\documentclass[11pt]{amsart}
\usepackage{amscd}
\usepackage{amsmath,amssymb,latexsym}
\usepackage{pb-diagram}
\usepackage{enumerate}
\usepackage{amsfonts}
\usepackage{dynkin-diagrams}

\newtheorem{thm}{Theorem} [section]
\newtheorem{Thm}[thm]{Theorem}
\newtheorem{Cor}[thm]{Corollary}
\newtheorem{Lem}[thm]{Lemma}
\newtheorem{Prop}[thm]{Proposition}

\newtheorem{Claim}[thm]{Claim}

\theoremstyle{definition}
\newtheorem{Def}[thm]{Definition}

\theoremstyle{remark}
\newtheorem{Remark}[thm]{Remark}

\numberwithin{equation}{section}

\newcommand{\la}{\langle}
\newcommand{\ra}{\rangle}
\newcommand{\beq}{\begin{equation}}
\newcommand{\eeq}{\end{equation}}
\newcommand{\be}{\begin{enumerate}}
\newcommand{\ee}{\end{enumerate}}
\newcommand{\bv}{\mathbb V}
\newcommand{\bw}{\mathbb W}
\newcommand{\bO}{\mathbb O}
\newcommand{\bG}{\mathbb G}
\newcommand{\C}{\mathbb C}

\newcommand{\bT}{\mathbb T}

\newcommand{\bB}{\mathbb B}
\newcommand{\bN}{\mathbb N}
\newcommand{\mS}{\mathcal S}
\newcommand{\mF}{\mathcal F}
\newcommand{\mR}{\mathcal R}

\newcommand{\mO}{\mathcal O}

\newcommand{\mW}{\mathcal W}
\newcommand{\mf}{\mathfrak}

\newcommand\N{\operatorname{N}}
\newcommand\T{\operatorname{T}}

\newcommand\ind{\operatorname{ind}}
\newcommand\Ind{\operatorname{Ind}}

\newcommand\Hom{\operatorname{Hom}}

\theoremstyle{plain}

\title{Fourier Transform and the minimal representaiton of $E_7$}
\author{Wee Teck Gan and Nadya Gurevich}
\address{W.T. Gan: Department of Mathematics, National University of Singapore, 10 Lower Kent
Ridge Road Singapore 119076
}
\email{matgwt@nus.edu.sg}
\address{N. Gurevich: Department of Mathematics, Ben-Gurion University of the Negev, Be'er Sheva,  Israel 8410501}
\email{ngur@math.bgu.ac.il }
\begin{document}

\maketitle

\section{Introduction}\label{sec:intro}
Let $F$ be a $p$-adic field and $H$ be the group of $F$-points of 
a reductive algebraic group, or a central covering of such a group.
The minimal representation $\Pi$ of $H$, which exists and is unique 
for split simply-connected simply-laced groups, is a necessary ingredient for 
the theta correspondence for various dual pairs $(G_1,G_2)$ in $H$. The theta correspondence, which studies the restriction of $\Pi$ to $G_1\times G_2$,
is an efficient method for constructing functorial lifting between representations of $G_1$ and $G_2$,
both locally and globally. To study the theta correspondence,
it is useful  (if not necessary) to have concrete models of the minimal representation. 

The unitary minimal representation
is often realized on the Hilbert space $L^2(X)$ of square integrable functions on a space $X$ and its subspace $\mS(X)$ of smooth vectors is the smooth minimal
representation of $H$, also denoted by $\Pi$. Indeed, it is typical to have several such models for the representation
$\Pi$. For each model, there  is usually a maximal parabolic subgroup $Q=MN$,
with Levi subgroup $M$ and unipotent radical $N$,
whose action is given by explicit geometric formulas.
To complete the description of the action of $H$, it is necessary and sufficient  to determine the action of an additional element $s \notin Q$.
In particular, for an involution $s \notin Q$ normalizing $M$, one may expect a nice formula for $\Pi(s)$. 

For example, in the context of classical theta correspondence, a minimal representation is a Weil representation of the double cover  
$\widetilde{Sp}(W)$ of a symplectic group $Sp(W)$.  Given a polarization $W=W^+\oplus W^-$, with associated Siegel parabolic subgroup $Q = MN$ stabilizing $W^+$,
the corresponding Schrodinger model of $\Pi$  is realized on the space of Schwarz functions $\mS_c(W^-)$ and a representative of the longest  element of the Weyl group $s$ acts by the classical Fourier transform.

For a split simply-laced group $H$, the unitary minimal representation 
has been constructed by Kazhdan and Savin in \cite{KaSa}, using 
a Heisenberg model, i.e. a model where the action of the standard Heisenberg 
parabolic group $P=LV$ is explicit. For the $L^2$-Heisenberg model, Kazhdan and Savin also described the action of an additional element $s \notin P$, corresponding to the reflection in the highest root.   
Later, Magaard-Savin \cite{MagaardSavin} and Kazhdan-Polishchuk \cite{KazhdanPolishchuk} provided a detailed analysis of  the smooth model of this Heisenberg model.    

On the other hand,  Savin \cite{Sa2} constructed the analog of the Schrodinger model for the minimal representation, on which a maximal parabolic subgroup $Q=MN$ with abelian radical $N$ acts explicitly. 
The space of this Schrodinger model is a space of smooth functions on a cone $X$.
The cone $X$ is the minimal nonzero $M$-orbit in the Lie  algebra $\bar {\mf{n}}$
of the opposite unipotent radical $\bar N$. For some groups, such as those of type $D_n$, there are several choices for $Q$ and hence several Schrodinger models. For others, e.g. $E_8$, there are no such parabolic 
subgroups with abelian unipotent radical and hence a Schrodinger model does not exist.

However, what has been lacking for a while is an explicit description of an additional element $s \notin Q$ (normalizing $M$). 
Thus, our knowledge of the Schrodinger model of the minimal representation has  for some time been less complete than desired, compared to the classical case of the Weil representation. In particular, one would like to have the analog of the Fourier transform on the cone $X$.

Several subsequent works have partially rectified this situation.
As an  example, when  $H=O_{2n}$,  $X$ is  the cone of isotropic vectors  in a quadratic space of dimension $2n-2$. The parabolic subgroup $Q$ with Levi subgroup $M= GL_1\times O_{2n-2}$ acts explicitly. In a recent paper \cite{GurevichKazhdan} of the second author with Kazhdan, the formula for the action  of an involutive element $s$ normalizing $Q$ and commuting with $O_{2n-2}$ was obtained. The resulting integral operator is the desired  Fourier transform on the cone $X$. We should mention that over archimedean fields, the analogous analysis of the Schrodinger model and  the formula for a  Fourier transform on a cone has been obtained in the works of Kobayashi-Mano \cite{KobayshiMano},  Hilgert-Kobayashi-M\"ollers-Orsted 
\cite{HKMO12}  and Hilgert-Kobayashi-M\"ollers 
\cite{HKM14}.

The goal of this paper is to complete the description of the Schrodinger model for the minimal representation of a split adjoint exceptional group $H$ of  type $E_7$. The $L^2$-model is  
 realized on a space of functions $L^2(\Omega)$, where $\Omega$
is a variety of rank-one elements in the exceptional Jordan algebra
$J= H_3(\bO)$. In this model, the parabolic subgroup $Q=MN$ with Levi factor
of type $E_6$ acts by explicit formulas; see (\ref{Schrodinger:action}).
We consider an involution $s$ which conjugates $Q$
to the opposite parabolic $\bar Q$ and define
an integral operator $\Phi$ on the space $\mS_c(\Omega)$ of smooth functions
of compact support. In  Theorem  \ref{main}, which is the main result of this paper, we show that the restriction of $\Pi(s)$ to $\mS_c(\Omega)$ equals $\Phi$.

In the last section, we relate the operator $\Pi(s)$ to a
generalized Fourier transform for the group $E_6$,
defined by Braverman and Kazhdan \cite{BravermanKazhdan02}.
This phenomena has already been
observed for the minimal representation of the double cover of the symplectic 
group and of the even orthogonal group, which explains the
similarity between the formula for the operator $\Phi$ and the Fourier
transform on a cone defined in \cite{GurevichKazhdan}.

Our main theorem has several applications. The space of smooth vectors
$\mS(\Omega)$ in $L^2(\Omega)$ is not easy to describe. 
We know that the space contains the space $\mS_c(\Omega)$ of 
compactly support functions and that each function in $\mS(\Omega)$
has bounded support. It is thus necessary to describe the space 
$[\mS(\Omega)]_0$ of germs of functions near $0$ (which is the vertex of the cone $\Omega$). This space affords 
an action of the Levi subgroup $M$ and is canonically isomorphic 
to Jacquet model $\Pi_N$. Using the formula for $\Pi(s)$
we describe the space $[\mS(\Omega)]_0$ 
as a direct sum of a one-dimensional representation and the minimal
representation of the Levi subgroup $M$ of type $E_6$.
The structure  of the Jacquet module, that was initially obtained in \cite{Sa2}
is an important ingredient in the computation of theta correspondence. 

The formula for $\Pi(s)$ has recently been  used in \cite{LeJun}.
Let  $X=Spin_9\backslash F_4$ be the exceptional spherical variety
of rank $1$. The associated group to the variety $X$ is $G_X=PGL_2$.
The $X$-distinguished representations of $F_4$ are expected to be
functorial lifts from generic representations of $PGL_2$.
This functorial lift is realized by the theta correspondence
for the dual pair $(PGL_2, F_4)$ in $H$. Under fixed embedding
$F_4$ is contained in $M$ and a representative of the non-trivial
Weyl element $w$ of $PGL_2$ is mapped to $s$. Let $N'$ denote the unipotent radical
of the Borel subgroup of $PGL_2$. 
The authors use the space $\mS_c(\Omega)$ to define a transfer map
$$t_\psi:\mS_c(\Omega_1)_{Spin_9}\rightarrow \mS((N',\psi)\backslash PGL_2)_{N',\psi},$$
satisfying a relative character identity. 
Here $\Omega_1$ is the set of elements in $\Omega$ of trace one
naturally identified with the space $Spin_9\backslash F_4$.
The formula for $\Pi(s)$ is essential for obtaining an explicit geometric formula for this transfer map.

We outline the content of the paper. In section \ref{sec:preparation},
we introduce notation for Jordan algebra and group structure of $H$.
In section \ref{sec:minimal}, we compare the Heisenberg and Schrodinger models of  the minimal
representation and write the transition map between them.
In section \ref{sec:s(omega)}, we prove that a space
of operators on $\mS_c(\Omega)$ with certain equivariant properties
is one-dimensional. 
In section \ref{sec:Phi}, we define a candidate formula for the
operator $\Pi(s)$ by an integral and prove its convergence.
In section \ref{sec:approximation}, we prove an auxiliary result
for the difference between integral of a function on a quadratic space
over a cone and a hyperboloid. 
In section \ref{sec:main} we prove the main theorem \ref{main}.
In section \ref{sec:jacquet}, we apply the formula for 
for the operator $\Phi$ to describe Jacquet module $\Pi_N$. Finally, in section \ref{sec:BK},
we relate the formula for $\Pi(s)$ to Braverman-Kazhdan operator
for a group $G$ of type $E_6$.

\section{Preparation}\label{sec:preparation}
In this preparatory section, we introduce the various notions and objects that will play a role in this paper: octonions, the exceptional Jordan algebra, the exceptional group $E_7$ and its Siegel and Heisenberg parabolic subgroups, along with other pertinent group theoretic structures. 

\subsection{Octonions}
The octonion algebra $\bO$ over a field $F$
is a (non-associative)  algebra 
of dimension $8$ with identity, equipped with 
a non-degenerate quadratic norm $\N:\bO\rightarrow F$, 
satisfying $\N(xy)=\N(x)\N(y)$.
The algebra admits a unique anti-involution $x\mapsto \bar x$
such that $\N(x)=x\bar x$. The trace $\T:\bO\rightarrow  F$
is defined by $\T(x)=x+\bar x$. The bilinear form defined by
$\la x,y\ra\rightarrow \T(x\bar y)$ is non-degenerate
and $\N(x)=\frac{1}{2}\la x,x\ra$. 
Let us record several properties of $\T$

\begin{enumerate}
\item
$\T(xy)=\T(yx).$
\item
$\T((xy)z)=\T(x(yz))$; in particular $\T(xyz)$ is well defined.
\end{enumerate}

\subsection{Albert algebra $J$}
Consider the cubic Jordan algebra $J=H_3(\bO)$ of $3\times 3$
Hermitian matrices over $\bO$. It consists of elements of the form
\begin{equation}\label{j:matrix}
j=j(a,b,c,x,y,z)=\left(\begin{array}{ccc}
  a & z & \bar y\\
\bar z & b & x\\
y& \bar x & c
\end{array}\right),\quad a,b,c\in F, x,y,z\in \bO.
\end{equation}
\begin{itemize}
    \item
The product defined by $A\circ B= \frac{1}{2}(AB+BA)$ is commutative.
\item
The identity matrix $I$ serves as the identity element. 
\item
There is a quadratic adjoint map $\#: J\rightarrow J$ defined by

\beq
\left(\begin{array}{ccc}
  a & z & \bar y\\
\bar z & b & x\\
y& \bar x & c
\end{array}\right)
\mapsto
\left(\begin{array}{ccc}
  bc-\N(x) & \overline{xy}-cz & zx-b\bar y\\
xy-c\bar z & ca-\N(y) &   \overline{yz}-ax\\
\overline{zx}-by& yz-a\bar x  & ab-\N(z)
\end{array}\right)
\end{equation}
\item
The cubic form $\bN:J\rightarrow F$ is defined as the determinant. Precisely
$$\bN(j(a,b,c,x,y,z))=abc-a\N(x)-b\N(y)-c\N(z)+\T(xyz).$$
\item
The trace $\bT:J\rightarrow F$ is defined by $\bT(j(a,b,c,x,y,z))=a+b+c.$
\item
There is a  natural pairing
$$\la\cdot,\cdot\ra_J:J\times J\rightarrow F,\quad
\la c_1,c_2\ra_J=\bT(c_1\circ c_2).$$

\item
Every element in $C\in J$ has a rank $0\le rank(C)\le 3$, defined as follows.
\begin{enumerate}
  \item If $C=0$, then $rank(C)=0$
  \item If  $C\neq 0, C^\#=0$, then $rank(C)=1$
    \item  If $C\neq 0, C^\#\neq 0$ and $\bN(C)=0$, then $rank(C)=2$
    \item If $\bN(C)\neq 0$, then $rank(C)=3$.
      \end{enumerate}
\end{itemize}
We denote by  $\Omega\subset J$ the set  of elements of rank $1$.
Its affine closure $\Omega^{cl}$ is the set of elements of rank at most $1$,
and $0$ is the unique boundary point. 
Let $\Omega^0\subset \Omega$ be the dense subset of elements
$j(a,b,c,x,y,z)\in \Omega$ such that $a\neq 0$.
Such element is completely determined
by three coordinates. We write explicitly

\begin{equation}\label{omega:matrix}
c(a,x,\bar y)=\left(\begin{array}{ccc}
       a &       x        &   \bar y\\
\bar x   & \frac{N(x)}{a} &    \frac{\overline{xy}}{a}\\
   y   &      \frac{xy}{a}         &      \frac{N(y)}{a}
\end{array}\right)\in \Omega^0,\quad a,\in F^\times , x,y\in \bO.
\end{equation}
The element $c_0=(1,0,0)\in \Omega^0$ will be often used in the paper. 

We shall further need an auxiliary set $\Omega^1$ defined by 
\beq \label{omega1}
\Omega^1=\{c(a,x,\bar y),\quad  a\in F^\times, x,y\in \bO^\times\}
\subset \Omega^0.
\eeq

The group $Aut(J,\bN)$ of automorphisms of
the vector space $J$, preserving
the norm $\bN$ is isomorphic to the split adjoint group of type $E_6$.
The action preserves the rank of the elements and acts transitively
on $\Omega$.

\subsection{The group $H$}
Let $H$ be the adjoint split group of type $E_7$ with Lie algebra
$\mf{h}$. We fix a maximal Cartan subalgebra $\mf{t}$,
maximal split torus $T$ and the absolute root system $R(H,T)$.
The choice of positive roots $R^+(H,T)$ in
$R(H,T)$ determines the Borel subgroup $B=T\cdot U$
and its opposite $\bar B=T\cdot \bar U$.
 There is a root space decomposition
 $$\mf{h}=\mf{t}\oplus( \oplus_{\alpha\in R(H,T)} \mf{h}_\alpha).$$ 

 For the labeling of simple roots
 \dynkin[labels={\alpha_1,\alpha_2,\alpha_3,\alpha_4,\alpha_5,\alpha_6,
     \alpha_7}]{E}{7}
we write $(n_1\ldots n_7)$ for the root
   $\sum_{\i=1}^7 n_i\alpha_i$. 
By $\omega_i, 1\le i\le 7$ we denote the fundamental weights.

For any root $\alpha\in R(H,T)$ we denote by $U_\alpha$ the corresponding
root subgroup. For a  subgroup $U'\subseteq U$ we write simply
$\alpha\in U'$, if $U_\alpha\subseteq U'$. 

Next we describe a relative root system.  
Consider a system of $3$ orthogonal roots, 
 $$\beta_1=(2234321),\quad \beta_2=(0112221), \quad \beta_3=(0000001).$$
The roots $\beta_1,\beta_2,\beta_3$ all belong to the
unipotent radical $N$, of the parabolic subgroup $Q$, defined by 
$\alpha_7$. 
Note that $\beta_1$ is the highest root in $R^+(H,T)$, 
$\beta_2$ is the highest root among all the roots in $R^+(H,T)$ that
orthogonal to $\beta_1$. The root $\beta_3$ is the highest root
in $N$ orthogonal to  $\beta_1$ and $\beta_2$. 

Fix  a $\mf{sl}_2$-triple $(e_i,f_i,h_i)$  
corresponding to the root $\beta_i$ for each $i$.
Let $\mf{t}_0=Span\{h_1,h_2,h_3\}\subset \mf{t}$
be a three-dimensional Cartan subalgebra
and $T_0\subset T$ be the corresponding torus in $H$.  Consider the
decomposition of $\mf{h}$ under the adjoint action of $\mf{t}_0$.
\beq 
\mf{h}=\mf{h}_0\bigoplus 
\left(\bigoplus_{\gamma\in R(H,T_0)} \mf{h}_\gamma\right).
\eeq
The derived algebra  $[\mf{h}_0,\mf{h}_0]$ has type $D_4$ generated by the roots
$\alpha_2,\alpha_3,\alpha_4,\alpha_5\in R(H,T)$. 
The relative  root system $R(H,T_0)$ is of type $C_3$.
Let $\gamma_1,\gamma_2,\gamma_3$ be simple relative roots with $\gamma_3$
the long one. We denote by $(n_1n_2n_3)$
the root $\sum_{i=1}^3 n_i \gamma_i\in R(H,T_0)$.
The root spaces corresponding to the long roots in $R^+(H,T_0)$
$$\mf{h}_{001}=\mf{h}_{\beta_3},\quad 
\mf{h}_{021}=\mf{h}_{\beta_2},\quad \mf{h}_{221}=\mf{h}_{\beta_1}$$
are one-dimensional.
The root spaces $\mf{h}_\gamma$ for short roots are all eight-dimensional.
For any $\gamma\in R(H,T_0)$ define
$$R_\gamma=\{\alpha\in R(H,T): \alpha|_{T^0}=\gamma\},
\quad U_\gamma=\Pi_{\alpha\in R_\gamma} U_\alpha.$$

\subsection{The Siegel parabolic subgroup }
Since the roots $\beta_i$ are orthogonal, the triple
\beq \label{efh:triple}
e=\sum e_i, \quad f=\sum f_i,\quad  h=\sum h_i
\eeq
is also a $\mf{sl}_2$ triple in $\mf{h}$.
The action of the element $ad(h)$ gives rise to a decomposition 
$\mf{h}=\mf{n}\oplus\mf{m}\oplus \overline{ \mf{n}}$
of eigenspaces with the  eigenvalues $2,0$ and $-2$ respectively.  
The algebra $\mf{m}\oplus \mf{n}$ is a maximal parabolic subalgebra
with $\mf{n}$ abelian. The associated parabolic subgroup
 $Q=M\cdot N$ is called the Siegel parabolic subgroup and it defined
by the root $\alpha_7$.
 We denote by $\bar Q=M\bar N$ the opposite
parabolic subgroup.

The nilpotent radical $\mf{n}$ admits a structure of Jordan algebra,
isomorphic to $J$, as described in \cite{KobayashiSavin}, Section $2$. 
We review the relevant structures, referring for proofs to \cite{KobayashiSavin}.

The product on $\mf{n}$ is defined by $$x\circ y=\frac{1}{2}[x,[f,y]].$$
The element $e$ is the neutral element
and the elements $e_1,e_2,e_3$ are orthogonal idempotents.
These idempotents define the Pierce decomposition
$$\mf{n}=\left(\oplus_{i=1}^3 \mf{n}_{ii}\right)\oplus
\left(\oplus_{1\le i\neq j\le 3} \mf{n}_{ij}\right),$$
where $\mf{n}_{ii}=Span(e_i)$ and
$$\mf{n}_{ij}=Span\{X_\alpha:
\la\alpha,\beta^\vee_i\ra=\la \alpha ,\beta^\vee_j\ra=1,
\la\alpha,\beta^\vee_k\ra=0, k\neq i,j\}.$$
It is easy to see that 
$$\mf{n}_{12}=\mf{h}_{121},\quad \mf{n}_{13}=\mf{h}_{111},\quad 
\mf{n}_{23}=\mf{h}_{011}.$$
A structure of a composition algebra $\bO$ is defined   
on each space $\mf{n}_{ij},$ $i\neq j$ using only the 
Lie brackets and the Killing form $\kappa$.

By Jacobson coordinatization theorem it follows that
$\mf{n}$ is isomorphic to $H_3(\bO)=J$. 
Each coordinate of a matrix in $H_3(\bO)$ corresponds to a relative
root space.  
The group $[M,M]$ is an adjoint group of type $E_6$
and conjugation by $g\in [M,M]$ preserves the norm $\bN$.
  Precisely,
  $$\bN(gxg^{-1})=\bN(x), \quad g\in [M,M], x\in \mf{n}=J.$$

\subsubsection{The pinning maps $e_\gamma$}
Let us define the  pinning maps $e_\gamma$ for the roots $\gamma\in R^+(H,T_0).$
We fix an isomorphism $J\mapsto \mf{n}$ 
such that each coordinate of a matrix in $J$ corresponds a relative
root space. In particular, for root $\gamma$ in $N$
there are  pinning maps 
$$\left\{\begin{array}{ll} e_\gamma:F\rightarrow U_\gamma & 
\gamma {\rm\, is\, long}\\
e_\gamma:\bO\rightarrow U_\gamma & 
\gamma {\rm\, is\, short}
\end{array}\right.
$$
such that 
$$j(a,b,c,x,y,z)\mapsto 
e_{221}(a)e_{021}(b)e_{001}(c)e_{011}(x)e_{111}(\bar y)e_{121}(z).$$
where $a,b,c\in F$ and $x,y,z\in \bO$.

Let us define the maps $e_\gamma:\bO\rightarrow U_\gamma$
for the roots $\gamma\in \{(100),(110),(110)\}$ such that $U_\gamma$ is not contained in  $N$.
Start with a short root $\gamma\in\{(100),(110)\}$. Note
that  $U_{\beta_1-\gamma}\subset  N$. 
There exist non-degenerate pairings
$$[\cdot, \cdot]:U_\gamma\times U_{\beta_1-\gamma}\rightarrow U_{\beta_1}\simeq F,
\quad \la \cdot,\cdot\ra:\bO\times \bO\rightarrow F, \quad \la x,y\ra= \T(x\bar y).$$
Hence the isomorphism $e_{\beta_1-\gamma}:\bO\rightarrow U_{\beta_1-\gamma}$
determines the  isomorphism $e_{\gamma}:\bO\rightarrow U_{\gamma}$
which preserves the pairing. Precisely
$$[e_{100}(x),e_{121}(y)]=e_{221}( \la x,y\ra)=[e_{110}(x),e_{111}(y)].$$
Similarly we define the isomorphism $e_{110}:\bO\rightarrow U_{110}$
such that 
$$[e_{010}(x),e_{011}(y)]=e_{021}( \la x,y\ra).$$

Finally, for any short  $\gamma\in R^+(H,T_0)$ we define
$e_{-\gamma}:\bO\rightarrow U_{-\gamma}$  such that
$\kappa(e_{-\gamma}(x),e_\gamma(y))=\T(xy)$ for all $x,y\in \bO$.

For any short root $\gamma$ we have $U_\gamma=\Pi_{\alpha\in R_\gamma} U_\alpha$.
This defines a decomposition $\bO=\oplus_{\alpha\in R_\gamma} L_\alpha$
into $8$ one dimensional spaces, such that $e_\gamma:L_\alpha\rightarrow U_\alpha$
is an isomorphism. Choosing a basis element $v_\alpha$ spanning $L_\alpha$ we get a pinning $e_\alpha:F\rightarrow U_\alpha$.
Since $H$ is simply-laced, whenever $\alpha,\beta,\alpha+\beta$ 
are positive roots in $R(H,T)$ one has 
$[e_\alpha(r),e_\beta(s)]=e_{\alpha+\beta}(N_{\alpha,\beta}\cdot rs)$ for some non-zero structure constants $N_{\alpha,\beta}$.

\subsection{Commutator relations}

In the following, we write several commutator relations in $H$. 
We begin with a simple lemma that holds for a general group. 

\begin{Lem}\label{comm:relation:abc}For any  elements $a,b,c$ in a group $G$
  one has 
  \begin{enumerate}
  \item
    $[ab,c]=[a,[b,c]]\cdot [b,c]\cdot [a,c],$
  \item
    $[c,ab]= [c,a]\cdot[c,b]\cdot [[b,c],a].$
  \end{enumerate}
  \end{Lem}
\begin{proof}
The first identity is obviously equivalent to the second. For the second
\begin{multline*} 
  [c,ab]=cabc^{-1}b^{-1}a^{-1}=cac^{-1}[c,b]a^{-1}=\\
  [c,a] a [c,b] a^{-1}=[c,a]\cdot [c,b]\cdot[[c,b]^{-1},a]=
  [c,a]\cdot[c,b]\cdot [[b,c],a].
    \end{multline*}
 \end{proof}

\begin{Prop}\label{pinning:identities}
  Let $s,r,z\in \bO$ and $c\in F$.
  \begin{enumerate}
    \item
  $[e_{100}(s),e_{021}(c)]=e_{121}(cs) e_{221}(N(s)c)$,

    \item
 $[e_{110}(r),e_{001}(b)]=e_{111}(br)e_{221}(N(r)b),$
\item $[e_{100}(r),e_{011}(z)]\in U_{111},$
\item $[e_{110}(s),e_{011}(z)]=e_{121}(s\bar z)\in U_{121},$
\item
  $[e_{100}(r)e_{110}(s), e_{011}(z)]\in U_{111}U_{121}\cdot
  e_{221}(\T(\bar r s \bar z))$.  
  \end{enumerate}
  \end{Prop}

\begin{proof}

For any $r\in \bO$ we write  $e_{100}(r)=\Pi_{\alpha\in R_{100}} e_\alpha(r_\alpha)$
for some $r_\alpha\in F$. 
Consider a partition $R_{100}=R^1\cup R^2$ with $|R^1|=|R^2|=4$,
defined by
$$R^1=\{\alpha\in R_{100}: n_2=0\}, \quad R^2=\{\alpha\in R_{100}: n_2=1\}.$$
This defines a partition $R_{121}=S^1\cup S^2$, where $S^i=R^i+(021)$ for $i=1,2$.

For any root $\alpha\in R^1$ there exists unique root 
$\alpha'\in R^2$ such that $\alpha+\alpha'+(021)=(221).$
For any two roots $\alpha,\alpha'\in R^i$, $\alpha+\alpha'+(021)$ is not a root.
It follows from Lemma \ref{comm:relation:abc}
that 

$$[\prod_{\alpha\in R^i} e_\alpha(r_\alpha), e_{021}(c)]=
\prod_{\alpha\in R^i} [e_\alpha(r_\alpha),e_{021}(c)]=\prod_{\alpha\in S^i} 
e_\alpha( N_\alpha r_\alpha c)$$
for some non-zero constants $N_\alpha$. 

We can write $e_{100}(r)=e_{100}(r_1)e_{100}(r_2),$
where $e_{100}(r_i)=\Pi_{\alpha\in R^i} e_\alpha(r_\alpha)$.

$$[e_{100}(r_1),[e_{100}(r_2),e_{021}(c)]]=
e_{221}
(\sum_{\alpha \in R^1} M_{\alpha,\alpha'}r_\alpha r_{\alpha'} c),$$
where $\alpha+\alpha'+(021)=(221)$ and $M_{\alpha,\alpha'}$ are non-zero constants. By Lemma \ref{comm:relation:abc} again we obtain
  $$[e_{100}(r),e_{021}(c)]=e_{121}(L(r)c) e_{221}(Q(r)c)),$$
  where $L:\bO\rightarrow \bO$ is a linear isomorphism
  and $Q$ is a non-degenerate quadratic form on $\bO$.
  We have to show that $L(r)=r$ and $Q(r)=\N(r)$ for all $r$.

  We use the fact that conjugation by $e_{100}(r)\in [M,M]$
  preserves the norm $\bN$ on $\mf{n}\simeq J$. One has
  $$e_{100}(r)(e_{001}(1)e_{021}(c))e_{100}(-r)=
  e_{001}(1)e_{021}(c)\cdot e_{121}(L(r)c)e_{221}(Q(r)c).$$
  This implies that $L(r)\overline{L(r)}=Q(r)$ for all $r$.

  We apply Part $(1)$ of Lemma \ref{comm:relation:abc}
  with $a=e_{100}(r), b=e_{100}(s),c=e_{021}(1)$. 
  $$[e_{100}(r+s),e_{021}(1)]=
  [e_{100}(r),[e_{100}(s),e_{021}(1)]]\cdot[e_{100}(r),e_{021}(1)]\cdot
  [e_{100}(s),e_{021}(1)].$$
  Comparing the argument  of $e_{221}$ for both sides we obtain
  $$Q(r+s)=Q(r)+Q(s)+\T(r\cdot \overline{L(s)})$$
  Since $Q(r)=\N(L(r))$ it follows  that
  $\T(r\cdot \overline{(L(s)-s)})=0$ for all $s,r$. Equivalently
  $L(r)=r$ and $Q(r)=\N(r)$ as required. This proves part $(1)$.
  The part $(2)$ is proved similarly. Part $3$ follows since 
  $[U_{111},U_{100}]=[U_{111},U_{011}]=1$. 

  To prove part $(4)$, note that for any $\alpha\in R_{110}$ there exists a unique $\alpha'\in R_{011}$
such that $\alpha+\alpha'\in R_{121}$ and any root in $R_{121}$
is of that form. Hence, there exists a bilinear nondegenerate map $$ B:\bO\times \bO\rightarrow \bO, \quad (r,z)\mapsto B(r,z)$$
such that $[e_{110}(r),e_{011}(z)]=e_{121}(B(r,z))$. The goal is to show
that $B(r, z)=r\bar z$. Consider an element
$n=e_{021}(1)e_{011}(z)\in N=J$ of norm zero.
Hence also $\bN(e_{110}(r)n e_{110}(-r))=0$.
One has
$$e_{110}(r)ne_{110}(-r)=e_{001}(1)e_{110}(z)e_{111}(r)e_{121}(B(r,z))
e_{221}(\N(r)c).$$
%
Hence 
$$\T((r\bar z)\cdot B(r,z))-\N(r\bar z )-\N(B(r, z))=0$$
which is equivalent to $\N(r\bar z-B(r,z))=0$ for all $z,r\in \bO$.
Hence $B(r,z)=r\bar z$ as required. 

  Part $(5)$ is the consequence of Parts $(3)$ and $(4)$.
\begin{multline*}
  [e_{100}(r)e_{110}(s),e_{011}(z)]=
  [e_{100}(r),[e_{110}(s),e_{011}(z)]]\cdot U_{111} U_{121}=\\
  e_{221}(\T(\bar rs \bar z))\cdot U_{111} U_{121}.
  \end{multline*}

\end{proof}

\begin{Cor}\label{cor:identities} Let $v=e_{100}(r)e_{011}(s)\in V_1$ and
  $n=e_{001}(b)e_{021}(c)e_{011}(z)\in N$ with $r,s,z\in \bO,$
  and $b,c\in F$. Then
  $$[v,n]\in U_{111}\cdot U_{121}\cdot e_{221}(\N(r)c+ \N(s)b+\T(\bar rs\bar z)).$$
\end{Cor}

\begin{proof}
Since $[v,n]\in N$ for all $v,n$ as above it follows by
Lemma \ref{comm:relation:abc} that
$[v,\cdot] :N\rightarrow N$ is a homomorphism. Hence 
$$[v,n]=[v,e_{001}(b)]\cdot[v,e_{021}(c)]\cdot[v,e_{110}(z)].$$

Applying the parts  $(1),(2),(5)$ of the Proposition \ref{pinning:identities} we get 
$$[v,e_{001}(b)]=[e_{110}(s),e_{001}(b)]\in U_{111}e_{221}(N(s)b),$$
$$[v,e_{021}(c)]=[e_{100}(r),e_{021}(c)]\in U_{111}e_{221}(N(r)c),$$
$$[v,e_{110}(z)]\in U_{111}U_{121}\cdot e_{221}(T(\bar rs \bar z)),$$
as required. 

\end{proof}

\subsection{A Heuristic Model}
In this subsection, we describe a heuristic matrix model for thinking about 
the Lie algebra of type $E_7$, its  various root subspaces and parabolic subalgebras 
we have introduced above.  Though non-rigorous, this gives a more visual presentation of the above discussion and is quite a good approximation to reality.
 
The simple Lie algebra $\mf{sp}_6$ of type $C_3$ can be described as the matrix algebra 
$$\mf{sp}_6=\left\{ A\in M_{6\times 6}(F), \quad JA+A^tJ=0\right\},$$
where $J=\left(\begin{smallmatrix}0 & J_3\\-J_3 & 0\end{smallmatrix}\right)$
and $J_3=\left(\begin{smallmatrix} 0& 0 &1\\0 & 1& 0\\ 1&0&0\end{smallmatrix}\right)$.
The Borel subalgebra consists of upper triangular matrices.
Each coordinate $(i,j)$ of the matrix with $i\neq j$ 
corresponds to a root space. The pinning maps can be easily written.
For example  the elements $e_{100}(r)e_{110}(s)$ and $e_{221}(a)e_{021}(b)e_{001}(c)e_{121}(z)e_{111}(y)e_{011}(x)$
are represented by the matrices 

$$ \left(\begin{smallmatrix}
0&r&s &  &&\\
0&0& 0& &&\\
0&0&0 & &&\\
&& & 0&0&-s\\
&& & 0&0&- r\\
&& & 0&0&0
\end{smallmatrix}\right),\quad  
\left(\begin{smallmatrix}
0&0&0 &  y &z &a\\
0&0& 0& x&b& z\\
0&0&0 & c& x& y\\
&& & 0&0&0\\
&& & 0&0&0\\
&& & 0&0&0
\end{smallmatrix}\right).$$

The missing entries in the matrices above are  all zero. 
This allows to compute easily the structure constants 
using the multiplication of matrices. 

We can write informally the elements in the  algebra $\mf{h}$ 
as $6\times 6$ matrices as above. Each entry corresponds to a relative 
root space. For long (respectively short) relative roots, the entry belongs  to the field $F$ (respectively the octonion algebra $\bO$). 
For example,  the elements $e_{100}(r)e_{110}(s)$ and $e_{221}(a)e_{021}(b)e_{001}(c)e_{121}(z)e_{111}(\bar y)e_{011}(x)$ with $r,s,x,y,z\in \bO$ and $a,b,c\in F$
are represented by matrices 

$$ \left(\begin{smallmatrix}
0&r&s &  &&\\
0&0& 0& &&\\
0&0&0 & &&\\
&& & 0&0&-\bar s\\
&& & 0&0&-\bar r\\
&& & 0&0&0
\end{smallmatrix}\right),\quad 
\left(\begin{smallmatrix}
0&0&0 & \bar y &z &a\\
0&0& 0& x&b&\bar z\\
0&0&0 & c&\bar x& y\\
&& & 0&0&0\\
&& & 0&0&0\\
&& & 0&0&0
\end{smallmatrix}\right).$$

Thus the subalgebra $\mf{n}$ is clearly isomorphic 
to the algebra $J$ of Hermitian $3\times 3$ matrices over $\bO$. 

The exponential map on matrices defines  maps 
$\mf{h}_\gamma\rightarrow U_\gamma\subset H$
for the relative roots $\gamma$. 
Remarkably, the commutator relations in Proposition \ref{pinning:identities} are reflected by matrix multiplication.
For example it is easy to see that for all $s\in \bO, c\in F$
$$ \left[\left(\begin{smallmatrix}
1&s&0 &  &&\\
0&1& 0& &&\\
0&0&1 & &&\\
&& & 1&0&0\\
&& & 0&1&-\bar s\\
&& & 0&0&1
\end{smallmatrix}\right),
\left(\begin{smallmatrix}
1&0&0 & 0 &0 &0\\
0&1& 0& 0&c&0\\
0&0&1 & 0& 0& 0\\
&& & 1&0&0\\
&& & 0&1&0\\
&& & 0&0&1
\end{smallmatrix}\right)\right]=\left(\begin{smallmatrix}
1&0&0 & 0 &sc &N(s)c\\
0&1& 0& 0&0&\overline{sc}\\ 
0&0&1 & 0&0&0\\
&& & 1&0&0\\
&& & 0&1&0\\
&& & 0&0&1
\end{smallmatrix}\right),
$$
that is exactly the commutator identity in Proposition \ref{pinning:identities}, Part $(1)$.

The maximal parabolic subgroup $Q=MN$ in $H$ consists of the matrices $\left(\begin{smallmatrix} A & B \\ 0 & A^\ast \end{smallmatrix}\right)$, that visually look like 
matrices in the standard Siegel parabolic subgroup of 
$PGSp_6$, the adjoint split group of type $C_3$. 
The same is true about the Heisenberg parabolic subgroups
in $H$ and $PGSp_6$.

\subsubsection{The Levi subgroup}
The Levi subgroup $M$ has a one-dimensional
center $T_1$ generated by the cocharacter $\lambda=\omega_7^\vee$.
Note that $Ad(\lambda(t))x=tx$ for $x\in \mf{n}$.

The adjoint action of $M$ on $\mf{n}$ and $\overline{\mf{n}}$
are irreducible representations
with highest weights $\omega_1$ and $\omega_6$ respectively,
that are dual to each other.
The action of $M$ on $\mf{n}=J$ and $\overline{\mf{n}}$
preserves the norm $\bN$
up to the similitude factor $\mu=2\omega_7$ and $\mu^{-1}$ respectively.

\subsubsection{The element $s$}
Let $\varphi:SL_2\rightarrow H$ be the map associated to the
$\mf{sl}_2$-triple $(e,f,h)$, defined in (\ref{efh:triple}) and let
$$s=\varphi\left(\begin{smallmatrix} 0 & 1\\ -1 & 0\end{smallmatrix}\right).$$
   Then $s$ is an involution in $H$, normalizing $M$ and acting
   on it by an outer automorphism. Moreover, $Ad(s)(\mf{n})=\overline{\mf{n}} $.

   We use $s$ to define an inner product on $\overline{\mf{n}}$ by

\beq \label{inner:prod:J:def}
\la x, y\ra_{\overline{\mf{n}}}=-\kappa(Ad(s)(x),y).
\eeq
   \begin{Lem}\label{pairing:trace}
     For all $x,y\in \mf{n}=J$ holds 
   $$\la x, y\ra_{\overline{\mf{n}}}=\bT(x\circ y).$$
\end{Lem}

   \begin{proof}
    The action of the element $s$ on $\mf{n}$ is given explicitly by 
     $$Ad(s)(x)=\frac{1}{2}[f,[f,x]].$$
By commutativity of the product in $J$ we have  
$$x\circ y= \frac{1}{2}[y,[f,x]].$$
Finally $\bT(z)=\kappa(f,z)$ for any $z\in \mf{n}$. Hence 
\begin{multline}
-\kappa(Ad(s)x,y)=-\frac{1}{2}\kappa([f,[f,x]],y)=\\
\frac{1}{2}\kappa(f, [y,[f,x]])=\frac{1}{2}\bT([y,[f,x]])=\bT(x\circ y).
\end{multline}
     \end{proof}
   
Since $\kappa$ is $H$-invariant form it follows that 
\beq  \label{inner:product:equi}
\la Ad(g) x, Ad(g^s) y\ra_{\overline{\mf{n}}}=
\la x,  y\ra_{\overline{\mf{n}}}, \quad x,y\in\overline{\mf{n}},
\eeq
where $g\in M$ and $g^s=sgs^{-1}$. 

\subsection{Heisenberg parabolic subgroup}
Let $P=L\cdot V$ be the Heisenberg parabolic subgroup of $H$,
defined by the root $\alpha_1$.
The group $V$ is a Heisenberg group with the one-dimensional
center $Z$ generated by the highest root $\beta_1$.

We define abelian subgroups $V_1,V_2\subset V$ where
$V_1= V\cap M$, $V_2Z=V\cap N$.
One has $V=V_1\cdot V_2\cdot Z$ and $[V_1,V_2]=Z$.
Identifying $V_1,V_2$ with their Lie algebras and
$Z$ with $F$ via the pinning map $e_{\beta_1}$,
the commutator on $V$ defines a symplectic form on $V/Z$,
so that $V_1,V_2$ are complementary Lagrangian subspaces.

The group $M\cap L$ acts on $V\cap N=V_2Z$, preserving $V_2$ and $Z$.
One has $M\cap L= T_1\times M_1$, where $T_1$ is the center of $M$
and $M_1=Stab_{M\cap L}(e_{\beta_1}(1))$.

\subsection{Correspondence of characters}
Let $\psi$ be a nontrivial additive character of $Z=U_{\beta_1}$,
such that the character $\psi\circ e_{\beta_1}$ of $F$,
which by abuse of notation is also denoted by $\psi$,
has conductor $\mO_F$, the ring of integers of $F$.
For any $a\in F^\times$, we define the character $\psi_a$
of $F$ by $\psi_a(x)=\psi(ax)$. 

Denote by $\mR$ the group of characters on the abelian group
$V\cap N=V_2Z$, whose restriction to $Z$ is nontrivial.
The elements  of $\mR$ are parameterized by  $V_1T_1$, as follows. Any element in $T_1$ has a form $\lambda(a), a\in F^\times .$
For any $(v,\lambda(a))\in V_1T_1$ the character $\Psi_{v,\lambda(a)}$ is defined 
$$\Psi_{v,\lambda(a)}(v_2z)=\psi_a([v,v_2]z)=
\psi(\lambda(a)[v,v_2]z\lambda(a)^{-1}).$$ 

We shall identify the abelian groups $N,\bar N$ with their
Lie algebras $\mf{n},\bar{\mf{n}}$ via the exponential map. 
The Pontryagin dual of $N$ is naturally identified with $\bar N$
using the Killing form. Hence any complex character
of $N=\mf{n}$ has form $\Psi_y(x)=\psi(\la y,x\ra)$ for some
$y\in \bar{\mf{n}}$. 
The characters associated to $y\in \Omega\subset \bar{\mf{n}}$
are called rank one characters. Clearly $y\in \Omega^0$ parameterize
rank-one characters whose restriction to $Z$ is nontrivial.

\begin{Def} The  map $j:\Omega^0\rightarrow V_1T_1$ is defined by
  $$j(y)=(v,t) \quad \Leftrightarrow\quad \Psi_y|_{V\cap N}=\Psi_{v,t}$$
 Using the heuristic model, or Proposition \ref{pinning:identities} for the rigorous proof it is easy to 
 write the map $j$ in coordinates. 
 $$j(c(a,x,\bar y))=e_{100}\left(\frac{x}{a}\right)
 e_{110}\left(\frac{\bar y}{a}\right)\cdot  \lambda(a).$$
In particular the map $j$ is  a bijection.
\end{Def}

The equivariance properties of the map $j$ are written below.

\begin{Lem}\label{equiv:j}
  For any $y\in \Omega$ with $j(y)=(v,t)\in V_1T_1$ one has 
  \begin{enumerate}
 \item If $v_1\in V_1$, then  $j(v_1^{-1}yv_1)=(v_1v,t),$
  \item
 If $m=m_1t_1\in M_1 \times T_1$, then 
  $j(m^{-1}ym)=(m_1^{-1}vm_1, tt_1),$
\item If $n\in N\cap L$, then  
  $t[v,n]t^{-1}\in V_2\cdot e_{221}(\la y, n\ra).$
  \end{enumerate}
  \end{Lem}

\begin{proof}
  \begin{enumerate}
\item   
\begin{multline*}
\Psi_{v_1^{-1}yv_1}(v_2z)=\psi(\la y, v_1v_2z v_1^{-1}\ra)=
\psi(\la y, v_2z[v_1,v_2]\ra)=\\
\psi(t[v,v_2] [v_1,v_2]zt^{-1})=\psi(t[v_1v,v_2]zt^{-1})=\Psi_{v_1v,t}(v_2z).
\end{multline*}
Hence $j(v_1^{-1}yv_1)=(v_1v,t).$
    \item
\begin{multline*}
      \Psi_{m^{-1}ym}(v_2z)=\psi(\la y, m v_2z m^{-1}\ra)=
      \Psi_y(mv_2zm^{-1})=\\
      \psi(t [v,m v_2m^{-1}]zt^{-1})=
\psi(tt_1[v,m_1v_2m_1^{-1}z] (tt_1)^{-1})=\Psi_{m_1^{-1}vm_1,tt_1}(v_2z).
\end{multline*}
Hence $j(m^{-1}ym)=(m_1^{-1}vm_1,tt_1)$. 
  
\item We compute independently both sides of the equality. 
 Denote $n=e_{001}(c)e_{021}(b)e_{011}(z)$ with $b,c\in F$ and $z\in \bO$.

 By explicit form one has 
  $$\la c(a,x,\bar y), n\ra=\frac{\N(y)c}{a}+\frac{\N(x)b}{a}+\frac{T(yxz)}{a}.$$

By definition of $j$ one has 
$$j(c(a,x,\bar y))=(e_{100}(\frac{x}{a})e_{110}(\frac{\bar y}{a}),\lambda(a)).$$

By Corollary \ref{cor:identities}

$$t[v,n]t^{-1}=a\cdot \left(\frac{\N(x)b}{a^2}+\frac{\N(y)c}{a^2}+
\frac{\T(\bar x\bar y \bar z)}{a^2}\right).$$

The claim follows since
$\T(\bar x\bar y \bar z)=\T(zyx)=\T(yxz)$. 

  \end{enumerate}

\end{proof}

\subsection{Measures}
The derived group $M^{der}$ of $M$  acts transitively on $\Omega$ and the stabilizer of an element 
is isomorphic to the derived group  $Q_1^{der}$ of the maximal parabolic subgroup
$P\cap Q$. Since both groups are unimodular, there exists unique, up to multiplication by a scalar,
$M^{der}$-invariant measure $\eta$ on $\Omega$.
In fact $M$ acts on it by the character $|\omega_7(\cdot)|^{-8}$,
which ensures that
the action of $M$ on $\mS_c(\Omega)$ in (\ref{Schrodinger:action})  is unitary.

The measure defines by restriction a functional on $\mS_c(\Omega^0)$
which can be written more explicitly. For any $f\in \mS_c(\Omega^0)$ holds

$$\int\limits_{\Omega}f(z)\eta(z)=
\int\limits_{F^\times} \int\limits_{\bO}\int\limits_{\bO}
f(c(a,x,\bar y))|a|^{-4} d^\times a dx dy.$$

\section{Minimal Representations: two models}\label{sec:minimal}
The simply connected group of type $E_7$ admits a unique unitary minimal representation with trivial central character, that 
gives rise to the minimal representation $\Pi$ of the adjoint group $H$.
In this section we compare two different models of $\Pi$,
the Heisenberg model $\Pi_h$ on which the action of the Heisenberg parabolic
subgroup $P$ is explicit, and
the Schrodinger model $\Pi_s$
on which the action of the Siegel parabolic subgroup
$Q$ is explicit. The comparison will allow us to write an action
of the elements $\Pi_s(s_2)$, $\Pi_s(s_3)$  and eventually
$\Pi_s(s)$. The space of smooth vectors in each model defines a smooth
irreducible representation of $H$,
also denoted by $\Pi_h$ or $\Pi_s$ respectively. The smooth minimal representation is unramified and its Satake parameter is given by the semisimple element in the dual group
$H^\vee=E^{sc}_7$ associated to the subregular orbit. 
The Heisenberg model was constructed in \cite{KaSa}, see also
\cite{GanSavin} for more details. The Schrodinger model was
constructed in \cite{Sa2}, see also \cite{SavinWoodbury}. 

The main Theorem \ref{main} concerns the action of the element $s$ on $\Pi_s$.
The formula for this action will be derived from our knowledge of the model $\Pi_h$. Hence, the model $\Pi_h$ is only used as auxiliary model. After this section, we use only the model $\Pi_s$ and denote it just by $\Pi$. 

\subsection{The Heisenberg model}
The group $V_2Z$ is a maximal abelian subgroup of $V$.
We fix an additive character $\psi$ of $Z$ and extend
it trivially  across $V_2$ 
to the character $\Psi_{0,1}$ of $V\cap N=V_2Z$.

\subsubsection{The Weil representation $\rho_\psi$}
By Stone-von Neumann theorem there exists unique unitary representation
$\rho_\psi$ of $V$ with central character $\psi$. It is isomorphic
to the unitary induction of $\ind^{V}_{V\cap N} \Psi_{0,1}$
and can be realized on the space $L^2(V_1)$.
It can further be extended to the representation of $[L,L]\cdot V$
also denoted by $\rho_\psi$. Let us write the action of certain elements
\begin{enumerate}
\item  
  $$\rho_\psi(g)\phi(v)=\left\{
  \begin{array}{ll}
    |\omega_7|^4(g)|\phi(g^{-1}vg) & g\in M\cap [L,L]\\
    \Psi_{0,1}([v,g])\phi(v) & g\in N\cap L
  \end{array}\right.
  $$
\item
  $$\rho_\psi(s_2)\phi(e_{100}(x)e_{110}(y))=
  \int\limits_{\bO}\phi(e_{100}(x')e_{110}(y))\psi(\la x',x\ra) dx'$$
\item
 $$\rho_\psi(s_3)\phi(e_{100}(x)e_{110}(y))=
  \int\limits_{\bO}\phi(e_{100}(x)e_{110}(y'))\psi(\la y',y\ra) dx'$$ 
\end{enumerate}

\subsubsection{The representation $\Pi_h$ of $Q$}

It is easy to check that $L=[L,L]T_1$.  
The unitary completion of $\ind^{Q}_{[L,L]V} \rho_\psi$
is realized on $L^2(V_1T_1)$
with an explicit action of the elements in $P\cap Q$.
In fact, see \cite{KaSa}, this representation can be extended to an irreducible
unitary representation of $H$.

\begin{Thm} There exists a model of the unitary minimal representation
  on $L^2(V_1T_1)$
  where the elements of $P\cap Q$ act by the following formula. 
$$
\left\{\begin{array}{ll}
\Pi(v)\phi(v_1,t_1)=\phi(vv_1,t_1) & v\in V_1\\
\Pi(v_2z)\phi(v_1,t_1)=\psi(t_1[v_1,v_2]zt_1^{-1})\phi(v_1,t_1)& v_2z\in V_2Z \\
\Pi(m)\phi(v_1,t_1)=|\omega_7(m)|^4 \phi(m^{-1}v_1m,t_1) & m\in M_1\\
\Pi(t)\phi(v_1,t_1)=\phi(v_1,tt_1)& t\in T_1\\
\Pi(n)\phi(v_1,t_1)=\Psi_{0,1}([v_1,t_1nt_1^{-1}])
\phi(v_1,t_1) & n\in N\cap L
\end{array}
\right..
$$
\end{Thm}

The action of the element $s_2,s_3\in L$ is given by
\beq \label{heis:s2}
\Pi_h(s_2)f(e_{100}(x)e_{110}(y),\lambda(a))=|a|^4\int\limits_{\bO}f(e_{100}(x')e_{110}(y),\lambda(a))
\psi(\la x',ax\ra) dx'.
\eeq
\beq \label{heis:s3}
\Pi_h(s_3)f(e_{100}(x)e_{110}(y),\lambda(a))=|a|^4\int\limits_{\bO}f(e_{100}(x)e_{110}(y'),\lambda(a))\psi(\la y',ay\ra) dy'.\eeq

\subsection{The Schrodinger model}
The minimal representation admits a realization on the space $L^2(\Omega,\eta)$
that we denote by $\Pi_s$.  In this model the Siegel parabolic $Q=MN$ acts
explicitly as follows

\beq \label{Schrodinger:action}
\left\{\begin{array}{ll}
\Pi_s(x)\phi(y)=  \Psi_y(x)\phi(y) & x\in N,\\
\Pi_s(m)\phi(y)=|\omega_7(m)|^{4}\phi(m^{-1}ym)& m\in M \\
\end{array}
\right..
\eeq 

For the proof see \cite{SavinWoodbury}.

\subsection{Transition between models}
The  spaces $L^2(V_1T_1)$ and  $L^2(\Omega,\eta)$
are models of the  minimal unitary representation $\Pi$.
Recall that in Section \ref{sec:preparation} we have constructed a bijection $j:\Omega^0\rightarrow V_1T_1$
such that $\Psi_y=\Psi_{v,t}$ as characters on $V\cap N$, whenever $j(y)=(v,t)$.
We use the map $j$ to construct 
 a $H$-equivariant transition map
$$j^\ast:L^2(V_1T_1)\rightarrow L^2(\Omega,\eta).$$

For any $y\in \Omega^0$ we write  $j(y)=(v,\lambda(a))$. Let us define
$$j^\ast(f)(y)=|a|^{-6}f(j(y))=|a|^{-6}f(v,\lambda(a)). $$
Clearly, the map is an isomorphism of vector spaces.

\begin{Thm}
  The map $j^\ast$ is a unitary isomorphism in the one-dimensional space 
  $\Hom_H(L^2(V_1T_1), L^2(\Omega,\eta))$.
\end{Thm}

\begin{proof}
  First let us show that the isomorphism $j^\ast$ is unitary.
  $$\|j^\ast(f)\|^2=\int\limits_{\Omega} |j^\ast(f)(y)|^2\eta(y)=
  \int\limits_{F^\times}\int\limits_{V_2}
  |a|^{-12}|f(v,\lambda(a))|^2 |a|^{-4} d^\times a  dv=\|f\|^2.$$
Since the restriction of the unitary representation $\Pi$ to $P\cap Q$ is irreducible,
it is enough to show that $j^\ast\in \Hom_{P\cap Q}(L^2(V_1T_1), L^2(\Omega,\eta))$.

Let us check for any $g\in P\cap Q=T_1\cdot M_1\cdot (N\cap L)\cdot V$ and $y\in \Omega^0$ one has
$$j^\ast(\Pi_h(g)f)(y)=\Pi_s(g)j^{\ast}(f)(y).$$

Below we assume that  $j(y)=(v,\lambda(a))$.
 \begin{enumerate}
 \item $g=v_2z\in V_2Z=V\cap N$.
\begin{multline*}
   j^\ast(\Pi_h(v_2z)f)(y)=|a|^{-6}\Pi_h(v_2z)f(v,\lambda(a))=
   |a|^{-6}\psi_a([v,v_2]z)f(v,\lambda(a))=\\
   |a|^{-6}\Psi_y(v_2z) f(v,\lambda(a))=\Pi_s(v_2z)j^\ast(f)(y).
 \end{multline*}  
\item $g=v_1\in V_1$. By Lemma \ref{equiv:j}, part $1$
      \begin{multline*}
        j^\ast(\Pi_h(v_1)f)(y)=|a|^{-6}\Pi_h(v_1)f(v,\lambda(a))=\\
        |a|^{-6}
  f(v v_1,\lambda(a))=j^\ast(f)(v_1^{-1} y v_1)=\Pi_s(v_1)j^\ast(f)(y).
\end{multline*}
\item  $g=m\in M_1$. By Lemma \ref{equiv:j}, part $2$
\begin{multline*}
  j^\ast(\Pi_h(m)f)(y)=|a|^{-6} |\omega_7(m)|^4 f(m^{-1}v m,\lambda(a))=\\
 |\omega_7(m)|^4 j^\ast(f)(m^{-1}ym)=\Pi_s(m)j^\ast(f)(y).  
\end{multline*}
\item $g=t=\lambda(b)\in T_1$. By Lemma \ref{equiv:j}, part $2$
\begin{multline*}
  j^\ast(\Pi_h(t)f)(y)=|a|^{-6}f(v,\lambda(ab))=\\
  |b|^6|ab|^{-6}f(v, \lambda(ab))=
  |b|^6j^\ast(f)(by)=\Pi_s(t)j^\ast(f)(y)
  \end{multline*}
\item $g=n\in N_1$. By Lemma \ref{equiv:j}, part $3$
\begin{multline*} 
 j^\ast(\Pi_h(n)f)(y)=|a|^{-6}\Psi_{0,1}(a[v,n])f(v,\lambda(a))=\\
  \psi(\la y,n\ra)f(v,\lambda(a))=\Pi_s(n)j^\ast(f)(y).
\end{multline*}
\end{enumerate}
  
\end{proof}

This allows us to compute the action of the elements $\Pi_s(s_2), \Pi_s(s_3)$
\begin{Cor}\label{Schrodinger:action:s2:s3}
 $$\Pi(s_2)(f)(c(a,x,\bar y))=|a|^{-4} 
  \int\limits_{\bO} f(c(a,x',\bar y)) \psi\left( \frac{\la x,x'\ra}{a}\right) dx'.$$
  $$\Pi(s_3)(f)(c(a,x,\bar y))=|a|^{-4} 
  \int\limits_{\bO} f(c(a,x,\bar y')) \psi\left( \frac{\la y,y'\ra}{a}\right) dx'.$$
\end{Cor}

We also fix an element $s_0\in M$ of order $2$, that is an involution
and $s_0s_2s_0=s_1$.  In heuristic model the element 
$s_0$ is of the form  
$$s_0=\left(\begin{smallmatrix}S & 0\\0 & S^\ast\end{smallmatrix}\right), \quad S=\left(\begin{smallmatrix} 0& 1 &0\\1 & 0& 0\\ 0&0&1\end{smallmatrix}\right), \quad S^\ast=J_3 S J_3^{-1}$$
Its action is given by

\beq \label{action:s0}
\Pi(s_0) f(c(a,x,\bar y))=
f\left(c(\frac{\N(x)}{a},\bar x , \frac{\bar x \bar y}{a})\right).
\eeq

\section{The space $\mS(\Omega)$}\label{sec:s(omega)}
Let  $\mS^\infty(\Omega)$ and  $\mS_c(\Omega)$)
denote  the space of locally constant functions and
locally constant functions of compact support on $\Omega$ respectively.
The group $Q$ acts on $\mS^\infty(\Omega)$ by the formulas 
(\ref{Schrodinger:action}).
The space $\mS_c(\Omega)$ is preserved by this action, which is unitary with respect
to the inner product defined by the measure $\eta$ on $\Omega$.
Obviously, $\mS_c(\Omega)$ is contained in $L^2(\Omega,\eta)$
and is dense in it.
\begin{Def}
  An operator $T:\mS_c(\Omega)\rightarrow \mS^\infty(\Omega)$
  is called $(M,s)$ equivariant if
$T\circ \Pi(g)= \Pi(g^s)\circ T$ for all $g\in M$. 
  The space of $(M,s)$ equivariant operators is denoted by
  $\Hom_{(M,s)}(\mS_c(\Omega),\mS^\infty(\Omega))$. 
\end{Def}

Recall that the center $T_1$ of $M$,
acts on $\mS^\infty(\Omega)$ by
$$\lambda(r)f(y)=|r|^6\cdot f(r\cdot y).$$
Define the space $\mS(\Omega,z)$ for any complex number $z$
$$\mS^\infty\supset \mS(\Omega, z)=
\{f: \lambda(r) f(\omega)=|r|^z f(\omega)\},
\quad r\in F^\times.$$
This space is isomorphic to the  space of coinvariants
$\mS_c(\Omega)_{T_1,|\cdot|^z}$.

Let $Q_6$ be the maximal parabolic subgroup of $M$, defined by the simple root $\alpha_6$.
The group $M$ acts transitively on $\Omega$.
The stabilizer of $e_{-\beta_3}(1)$ is $U_6M_6,$
where $M_6\times T_1$  is the Levi subgroup of $Q_6$.
There is  an isomorphism of $M$-representations
$$\mS_c(\Omega)_{T_1,|\cdot|^z}\simeq
|\omega_7|^{-2(6-z)}\otimes \Ind^M_{Q_6} |\omega_6|^{-z}.$$


The following proposition will be used in the proof of Theorem \ref{main}.
\begin{Prop}\label{one:dim}
  The space  $\Hom_{(M,s)}(\mS_c(\Omega), \mS(\Omega,3))$
 is one-dimensional and the image of any non-zero operator
 in this space  is $|\omega_7|^{-6}\otimes \Pi_1$, where
 $\Pi_1|_{[M,M]}$ is the minimal representation of the semisimple group $[M,M]$
 of type $E_6$.
\end{Prop}  

\begin{proof}
  Any map $T\in\Hom_{(M,s)}(\mS_c(\Omega), \mS(\Omega,3))$ factors through
  the space of coinvariants $\mS_c(\Omega)_{T_1, |\cdot|^3}$ and hence
  gives rise to 
  $\tilde T\in \Hom_{(M,s)}(\mS(\Omega,-3),\mS(\Omega,3))$.
  The spaces $\mS(\Omega,3)$ and $\mS(\Omega,-3)$
  are identified with $|\omega_7|^{-6}\otimes \Ind^M_{Q_6} |\omega_6|^{-3}$
  and  $|\omega_7|^{-18}\otimes \Ind^M_{Q_6} |\omega_6|^{3}$ respectively.
\begin{Claim} $\Ind^M_{Q_6} |\omega_6|^{-3}$ is indecomposable 
of length $2$  with two distinct constituents and the restriction of the unique subrepresentation to $[M,M]$ is the minimal representation $\Pi_1$ of the adjoint group of type $E_6$.
\end{Claim}

The claim implies the proposition. 
  Indeed,  $Ad(s)(Q_6)=\bar Q_6$.
  and the map $Ad(s)$ is $(M,s)$-equivariant isomorphism
  $$\Ind^M_{\bar Q_6}|\omega_6|^{-3}\rightarrow \Ind^M_{Q_6}|\omega_6|^{3}$$
  Hence $$Ad(s)\circ \tilde T\in 
  \Hom_{M}(\Ind^M_{\bar Q_6}|\omega_6|^{-3},
  \Ind^M_{Q_6}|\omega_6|^{-3})$$
  This $\Hom$ space is at most one-dimensional by Claim. 
  It is exactly one-dimensional space,
  since the standard intertwining operator belongs to it and does not vanish on the spherical vector. 

  The claim follows from \cite{Weissman}
  for the group $[M,M]$ of type $E_6$. The representation $\Ind^{[M,M]}_{[M,M]\cap Q_6}|\omega_6|^{-3}$
is indecomposable  of length $2$
  and contains the minimal representation $\Pi_1$ of $[M,M]$ as unique subrepresentation by Corollary $6.3.2$ of $loc.cit$.
  
  Since $Q_6[M,M]=M$, the restriction of functions defines an isomorphism of the spaces
  $$|\omega_7|^{-6}\otimes\Ind^M_{Q_6}|\omega_6|^{-3}\rightarrow
 |\omega_7|^{-6}\otimes \Ind^{[M,M]}_{[M,M]\cap Q_6}|\omega_6|^{-3}$$
that is $[M,M]T_1$ equivariant. 
The representation $|\omega_7|^{-6}\otimes \Ind^M_{Q_1}|\omega_6|^{-3}$ of $M$
  is  of length at most  $2$. If it was irreducible, then
  its restriction to the normal group of finite index $[M,M]T_1$
  would be completely reducible, that is a contradiction.
  Hence it is of length two and the claim follows. 

  To sum up, the space $\Hom_{(M,s)}(\mS_c(\Omega), \mS(\Omega,3))$
  is at most one-dimensional and 
  any non-zero map must have image isomorphic to  $|\omega_7|^{-6}\otimes \Pi_1$

\end{proof}

Let us consider  the space  of smooth vectors $\mS(\Omega)$
in $L^2(\Omega,\eta)$,
which is the irreducible smooth minimal representation $\Pi$ of $H(F)$. 
It is shown in \cite{SavinWoodbury} that the space $\mS(\Omega)$ is contained in the
space of smooth functions of bounded support and  contains the space
$\mS_c(\Omega)$.
Clearly, the restriction of
$\Pi(s)$ to $\mS_c(\Omega)$ is a $(M,s)$-equivariant operator. 

\section{The operator $\Phi$}\label{sec:Phi}

The operator $\Pi(s):L^2(\Omega)\rightarrow L^2(\Omega)$ is
unitary $(M,s)$-equivariant of order $2$.
We start by writing a candidate formula for the
restriction of this operator to $\mS_c(\Omega)$, which we denote by $\Phi$.

The motivation for this formula comes from the study of Schrodinger
model for minimal representations for other groups. Specifically, for
the symplectic space with polarization $W=W_1\oplus W_2$
the minimal representation of the metaplectic group
$Mp_{2n}(W)$ is realized on the space $\mS_c(W_1)$
and  the analogous operator  $\Pi(s)$ acts by the classical Fourier transform.
For the orthogonal group $O(V)$ for a split quadratic space $V$
of dimension $2n+2$
the minimal representation is realized on $L^2$ functions on the
isotropic cone $X$ in the split $2n$-dimensional space. 
The operator $\Pi(s)$, studied in
\cite{GurevichKazhdan}, acts by a kernel operator, where the kernel distribution
is given by
$$K(x,y)=\int\limits_{F^\times} \psi(r\cdot \la x,y\ra+r^{-1}) |r|^{n-1} d^\times r,
\quad x,y\in X.$$

The formula for $\Pi(s)$ in our case is a slight variation
on the formula above, the power $n-1$ is replaced by $3$. 

In all the cases above the operator $\Pi(s)$ is closely related
to the generalized Fourier operators, defined by Braverman-Kazhdan.
We shall explain the relation in section \ref{sec:BK}.

\begin{Def} The operator $\Phi:\mS_c(\Omega)\rightarrow \mS^\infty(\Omega)$
  is defined by 
\begin{equation}\label{Phi:def}
  \Phi(f)(c_1)=
  \int\limits_{F^\times}
  \left(\int\limits_\Omega f(c_2)\psi(r\cdot \la c_1,c_2\ra)
  \eta(c_2)\right)\cdot 
  \psi(r^{-1}) |r|^3 d^\times r.
  \end{equation}
\end{Def}
Here $\la c_1,c_2\ra=\la c_1,c_2\ra_J.$ We will omit the subscript $J$ if there is no confusion. 
Since $f$ is of compact support, the inner integral,
denoted by $\hat \mR(f)(r c_1),$ converges absolutely.
The absolute convergence of
$$\int\limits_{F^\times} \hat \mR(f)(r c_1)\psi(r^{-1}) |r|^3 d^\times r$$
follows from the Proposition \ref{hat:R:properties} below.

\subsection{Convergence of the integral defining $\Phi(f)$}

We first prove several properties of the operator
$\hat \mR:\mS_c(\Omega)\rightarrow\mS^\infty(\Omega)$, defined by
$$\hat \mR(f)(c_1)=
\int\limits_\Omega f(c_2)\psi(\la c_1,c_2\ra)  \eta(c_2).$$
In particular, parts $(2),(3)$ of the next Proposition
imply the convergence of the integral in (\ref{Phi:def}).
\begin{Prop}\label{hat:R:properties}
  Let $f\in \mS_c(\Omega)$ and $c\in \Omega$ 
  \be
\item $\hat \mR$ is $(M,s)$ equivariant.  
\item $\hat \mR(f)$ is locally constant on the affine closure $\Omega^{cl}$ of $\Omega$.
  One has $\hat \mR(f)|_{V\backslash \{0\}}=\int\limits_{\Omega} f(c)\eta(c)$
  for an open neighborhood $V\subset \Omega^{cl}$ of $0$. 
\item The function on $F^\times$, defined by $r\mapsto \hat \mR(f)(r\cdot c)|r|^{4}$ is bounded.

  \ee 
\end{Prop}

\begin{proof}
  The first property follows from the facts that
  $$\la g c_1, g^sc_2\ra=\la c_1,c_2\ra,\quad 
|\omega_7(g^s)|=|\omega_7(g)|^{-1}, g\in M,$$
see Section \ref{sec:preparation}.

  Since $f$ is of compact support, there exists an open neighborhood $V$ 
  of $0$ such that $c_2\in supp(f)$ and $c_1\in V$ implies
   $\psi(\la c_1,c_2\ra)=1$.
  This implies part $(2)$.
 
It is enough to show part $(3)$ for $c=c_0=c(1,0,0)$. 
We shall compute it explicitly.
One has  $\la c_0,c(a,x,y)\ra=a$. Hence
  $$\hat \mR(f)(r\cdot c_0)=
  \int\limits_{F^\times}
  \left(\int\limits_{\bO}\int\limits_{\bO}
    f(c(a,x,y))|a|^{-5}dxdy\right)  \psi(ra)da.$$

    In particular $r\mapsto \hat \mR(f)(r c_0)$ is the one-dimensional
    Fourier transform of the function 
$$a\mapsto \mR(f)(a)=\int\limits_{\bO}\int\limits_{\bO}
      f(c(a,x,y))|a|^{-5}dxdy=|a|^3\Pi(s_2s_3)f(a c_0).$$
      
    The function $\mR(f)$ is clearly locally constant on $F^\times$ and
    vanishes for $|a|$ large. The following lemma describes
    it behavior for $|a|$ small.

    \begin{Lem} For any $f\in \mS_c(\Omega)$
      there exists constants $K_1,K_2$ and $\epsilon >0$ such that for
      $|a|<\epsilon$ such that  $\mR(f)(a)=K_1|a|^3+K_2.$
\end{Lem}

    \begin{proof}      
      Let $f=f^0$ be the normalized spherical function in $\mS_c(\Omega)$,
       computed in \cite{SavinWoodbury}. One has
      $f(c)=\zeta(-3)+\zeta(3)\|c\|^{-3}$ for $\|c\|\le 1$ and
      zero otherwise. Here $\|\cdot\|$ is the
      maximum norm on the vector space $\overline{\mf{n}}$.
      Let $c_0=c(1,0,0).$
      Hence
      $$\mR(f^0)(a)=|a|^3\Pi(s_2s_3)f^0(a c_0)=\zeta(-3)|a|^3+\zeta(3)$$
      for $|a|\le 1$.
      The vector $f^0$ generates $\mS(\Omega)$. By Iwasawa decomposition
      $H=N\cdot M\cdot H(\mO)$, 
     for any vector $f\in \mS_c(\Omega)$ one has
      $\Pi(s_2s_3)f=\sum_{i=1}^k \Pi(n_im_i)f^0$ for some $n_i\in N, m_i\in M$. 
      Hence
      $$|a^3|\Pi(s_2s_3)f(a c_0)=\sum_{i=1}^k |a|^3 \Psi_{n_i}(a c_0)
      |\omega(m_i)|^4f^0(a m_i^{-1} c_0 m_i).$$
      For $|a|$ small one has $\Psi_{n_i}(a c_0)=1$ for all $i$
      and $\|a m_i^{-1} c_0 m_i\|\le 1$. Hence it equals
      $K_1|a|^3+K_2$ for some constants $K_1,K_2$ depending
      on $f$. 
 \end{proof}
    For any $g\in \mS^\infty (F^\times)$
    such that  $g(a)=|a|^s$ for $|a|$ small and $g(a)=0$ for
    $|a|$ large, the function 
    $\hat g(a)|a|^{s+1}$ is bounded, where $\hat g$ is the 
    classical Fourier transform. 
    The proposition follows.

\end{proof}
  
The operator $\Phi$ is $(M,s)$ equivariant. This follows from
the property of equivariance of the measure $\eta$
and (\ref{inner:product:equi}).

Theorem \ref{main} implies that for any $f\in \mS_c(\Omega)$ one has
$\Phi(f)=\Pi(s)(f)\in \mS(\Omega)$ and hence has bounded support.
Let us explore the behavior of $\Phi(f)$ near zero. 

We consider a decomposition $\Phi=\Phi_1+\Phi_2$ where
\beq \label{Phi1}
\Phi_1(f)(c_1)=\int\limits_{F^\times}
  \left(\int\limits_\Omega f(c_2)\psi(\la c_1,c_2\ra \cdot r),
  \eta(c_2)\right)\cdot |r|^3 d^\times r,
  \eeq
\beq \label{Phi2}
  \Phi_2(f)(c_1)=\int\limits_{F^\times}
  \left(\int\limits_\Omega f(c_2)\psi(\la c_1,c_2\ra \cdot r)
  \eta(c_2)\right)\cdot (\psi(r^{-1})-1)|r|^3 d^\times r.
  \eeq

  Both operators are well-defined and are $(M,s)$ equivariant.
  Moreover the operator $\Phi_1$ belongs to the one-dimensional
  space $\Hom_{(M,s)}(\mS_c(\Omega), \mS(\Omega,3))$. This decomposition will be
  essential in the proof.  The next Lemma will also be used in the proof of Theorem \ref{main}.
  
\begin{Lem}\label{f0:lemma}
  Define a function $f_0\in \mS_c(\Omega)$ by
  $$f_0(w)=\left\{\begin{array}{ll}
  1 & w=c(a,x,\bar y), |a|=1, \|x\|, \|y\|\le 1\\
  0 & {\rm otherwise}
  \end{array}\right..$$
  Let $c_0=c(1,0,0)\in \Omega$.
  Then  $\Phi(f_0)(sc_0)=0$ for $|s|\gg 1$ and  $\Phi_1(f_0)(c_0)\neq 0$.
\end{Lem}

\begin{proof}
It follows by definition of $f_0$ that 
  $$\mR(a)(f_0)(c(1,0,0))=
  \int\limits_{F^\times}\int\limits_{\bO}\int\limits_{\bO}
  f_0(c(a,x,y)) |a|^{-4} d^\times a dx dy=\left\{
  \begin{array}{ll} 1 & |a|=1\\
    0 & |a|\neq 1
    \end{array}\right..
    $$

Hence 
$$\hat \mR(f_0)(rc(1,0,0))=\int\limits_{|a|=1} \psi(ra) da
=\left\{
\begin{array}{ll} 0 & |r|>q\\
-1 & |r|=q\\
1 & |r|\le 1
\end{array}\right. .
$$

$$\Phi(f_0)(sc(1,0,0))=\int\limits_F
\hat \mR(f_0)(sr\cdot c(1,0,0)) \psi(r^{-1})|r|^3 d^\times r=$$

$$\int\limits_{|r|\ge |s|} \psi(r)|r|^{-3} d^\times r -
\int\limits_{|r|= |s|q^{-1}} \psi(r)|r|^{-3} d^\times r,$$
which is zero for $|s|$ large.

On the other hand the value $\Phi_1(f)(c(1,0,0))$ equals
$$\int\limits_{F^\times}\hat \mR(f_0)(rc(1,0,0))|r|^3 d^\times r=
\int\limits_{|r|\le 1} |r|^3 d^\times r-
\int\limits_{|r|=q} |r|^3 d^\times r \neq 0,$$
as required.

\end{proof}

\begin{Remark}
The function  $\mR(\cdot )(f)$ is related to the Radon transform 
on $\Omega$. 
  For any $w\in \Omega$ consider a fibration $p_w:\Omega\rightarrow F$ 
  defined by $p_w(v)=\la w,v\ra$. The fibers $\Omega_{w,a}=\{v: \la v,w\ra=a\}$ are
  smooth for $a\in F^\times$ and admit measures $\eta_{w,a}$
  satisfying 
  $$I_\Omega(f):=\int\limits_{\Omega}f(x) \eta(x)=\int\limits_{F}
\int\limits_{\Omega_{w,a}}f(x) \eta_{w,a}(x) da, \quad f\in \mS_c(\Omega).$$
  
For any $a\in F^\times$ the Radon transform 
${\bf R}(a):\mS_c(\Omega)\rightarrow \mS^\infty(\Omega)$ 
is defined by 
$${\bf R}(a)(f)(w)=\int\limits_{\la v,w\ra=a} f(v)\eta_{w,a}(v).$$
It is easy to see that $\mR(a)(f)={\bf R}(a)(f)(c_0)$. The Radon 
transforms will also appear in Section \ref{sec:BK}.

\end{Remark}

\section{Approximation of a hyperboloid by a cone}\label{sec:approximation}

In this section we prove an identity between  distributions on quadratic
vector spaces that
is essentially used in the proof of Theorem \ref{main}. Although
the identity will be applied for the  quadratic space $(\bO,N)$
of dimension eight only, we work here with general split quadratic spaces. 

 Remind that $\psi$ is an additive character $F$ of conductor 
$\mO_F$ and  $\psi_a$ is the twisted character 
$\psi_a(x)=\psi(ax)$ for any $a\in F$. In particular 
$\psi_0$ is the trivial character of $F$. 

Let $(V,q)$ be a non-degenerate split
quadratic space of dimension $2n\ge 4$.
We denote by $V(s)$ the  fiber of $q$ over $s$. The fiber
$V(s)$ is a hyperboloid and is smooth
for $s\neq 0$. The cone $V(0)$ of isotropic vectors 
has  unique singular point $v=0$.
The differential forms $dv$, $dx$, such that the measures
$|dv|$ on $V$ and $|dx|$ on $F$  are self dual with respect to $\psi$,
give rise to differential forms $\eta_s$ on the smooth fibers, and 
to the measures $|\eta_s|$.  It is easy to show that
for $\dim(V)\ge 4,$
the integral over  the measure $|\eta_0|$ converges  for all
functions in $\mS_c(V)$. In the integrals below we will
write $\eta(x)$ and $dv$ instead of $|\eta(x)|$ and $|dv|$.

Consider a family of  distributions  $\delta_{s}$  on $\mS_c(V)$
given by $\delta_{s}(f)=\int\limits_{V(s)} f(v)\eta_s(v).$
By Fubini theorem for any $f\in \mS_c(V)$ holds
  $$\int\limits_F \delta_s(f) ds=\int\limits_{V} f(v) dv.$$ 
We denote by $\mF_V$ the Fourier transform with respect to
the character $\psi$, the form $q$  and the self-dual measure $dv$.

The integration over the fiber $V(s)$ can be  approximated by
the integration over the fiber $V(0)$. Precisely,

\begin{Thm}\label{delta:decomp}
  $$\delta_s(f)=\delta_0(f)+
  |s|^{n-1}\int\limits_{V} \mF_V(f)(v) H_{n-1}(sq(v))dv,$$
 where for $Re(z)>0$
  $$H_{z}(s)=\int\limits_{F^\times} \psi(st)(\psi(t^{-1})-1)|t|^{z} d^\times t.$$
\end{Thm}

The proof of this theorem uses the dual pair $SL_2\times O(V)$
inside $Sp_{2n}$.

For this section only we denote by $B=TN$ the Borel subgroup of $SL_2$. 
Let us introduce notations for some elements in $SL_2$.
$$w=\left(\begin{array}{ll} 0 & 1\\-1 & 0\end{array}\right), \quad
  n(s)=\left(\begin{array}{ll} 1 & s\\0 & 1\end{array}\right),\quad
h(a)= \left(\begin{array}{ll} a & 0\\0 & a^{-1}\end{array}\right),   $$

It is easy to check the relation
$$wn(-s^{-1})w= n(s)  h(s)\cdot  wn(s).$$

For any $z\in \C$ define the character $\chi_z$ of $T$ by
$\chi_z(h(a))=|a|^{z}$. 

Let  $\Ind_B(\chi_z)$ be the normalized principal series of $SL_2$
induced from the character $\chi_z$. In particular, for  any
$\varphi\in \Ind(\chi_z)$ one has 
\begin{equation}\label{inverse:sl2}
  (w\varphi)(wn(-s^{-1}))=|s|^{z+1} \varphi(wn(s))
  \end{equation}

For any $s\in F$ define
  the Whittaker functional $W_s\in \Hom_{N}(\Ind_B(\chi_z), \C_{\psi_s})$ by
$$W_s(\varphi)=\int_F \varphi(wn(t)) \psi(-st)dt$$
The integral $W_s(\varphi)$ converges absolutely for $Re(z)>1$ and stabilizes
for all $Re(z)>0$ for $s\neq 0$.
Conversely,
\begin{equation}\label{W:to:f}
  \varphi(wn(s))=\int_F W_r(\varphi)\psi(rs) dr
  \end{equation}

\begin{Lem}
 For any $\varphi$ in a principal series  $\Ind_B(\chi_z)$ one has 
 \begin{equation}\label{asym:sl2}
  W_s(\varphi)-W_0(\varphi)=
  |s|^z\int\limits_F W_r( w\varphi) H_z(sr) dr
  \end{equation}
\end{Lem}
\begin{proof}
On has 
$$W_s(\varphi)-W_0(\varphi)= \int\limits_F \varphi(wn(t))(\psi(-st)-1)dt,$$
which equals by \ref{inverse:sl2} to
$$\int\limits_F (w\varphi)(wn(-t^{-1}))|t|^{-z-1}(\psi(-st)-1)dt$$
We apply \ref{W:to:f}
  $$\int\limits_F \int\limits_{F} W_r(w\varphi)
  \psi(-rt^{-1})dr |t|^{-z}(\psi(-st)-1)d^\times t.$$
Making change of variables $t\mapsto -(st)^{-1}$ we obtain
 $$|s|^z \int\limits_F \int\limits_{F} W_r(w\varphi)
  \psi(rst)(\psi(t^{-1})-1)|t|^z dr  d^\times t =$$

  $$|s|^z \int\limits_F W_r(w\varphi)\left( \int\limits_{F} 
  \psi(rst)(\psi(t^{-1})-1)|t|^z dr  d^\times t\right)  dr$$
  The inner integral converges absolutely and equals $H_z(rs)$.
  Hence the above equals 
  $$|s|^z \int\limits_F W_r(w\varphi)H_z(rs) dr$$
  as required. 
   
\end{proof}
  
 The group  $SL_2\times O(V,q)$ acts on $\mS_c(V)$ via the Weil
 representation $\omega_{\psi,q}$. We shall use the following formula

 $$\omega_{\psi,q}(w) f(v)=\mF_V(f),  \quad \omega_{\psi,q}(n(s))f(v)=
 \psi(sq(v))f(v).$$
 
  The Rallis map
  $$R: S_c(V)\rightarrow \Ind_B(\chi_{n-1}),
  \quad R(f)(g)=\omega_{\psi,q}(g)f(0)$$
  is $O(V,q)$-equivariant.

\begin{Lem}
  One has \begin{equation}\label{delta:Ws}
    \delta_s(f)=W_s(R(f)).\end{equation}
\end{Lem}

\begin{proof}
  One has
  $$R(f)(wn(s))=\omega_{\psi,q}(wn(s))f(0)=
  \int\limits_V f(v)\psi(s q(v)) dv=
  \int_F \delta_t(f)\psi(ts) dt.$$
Applying the inverse Fourier transform we get $\delta_s(f)=W_s(R(f)).$
\end{proof}

\begin{proof} of Theorem. 
Applying (\ref{delta:Ws}) and (\ref{asym:sl2}) for $\varphi=R(f)$ we get
$$\delta_s(f)-\delta_0(f)=|s|^{n-1}\int\limits W_r(wR(f))H_{n-1}(rs) dr
=|s|^{n-1}\int\limits W_r(R(wf))H_{n-1}(rs) dr=$$
$$|s|^{n-1}\int\limits_F \delta_r(\mF(f))H_{n-1}(rs)dr =
|s|^{n-1}\int\limits_V \mF_V(f)(v)\,  H_{n-1}(q(v)s) dv$$
as required.
\end{proof}

Theorem \ref{delta:decomp} will be applied to the quadratic space
$(\bO,N)$, so that $n=4$.
Precisely, denoting by $\bO(s)$ the octonions  of norm $s$ we obtain 

\begin{equation}\label{O(x):decomp}
  \int\limits_{\bO(s)}f (x)\eta_s(x)=\int\limits_{\bO(0)}f (x)\eta_0(x)+
  |s|^{n-1}\int\limits_{V} \mF_\bO(f)(v) H_{3}(sq(v))dv,
  \end{equation}
 where 
  $$H_{3}(s)=\int\limits_{F^\times} \psi(st)(\psi(t^{-1})-1)|t|^{3} d^\times t.$$
Let us mention another simple identity of functions
on $\bO$, that will be used later.

\begin{Lem}\label{Fourier:O:shift}
  Let $g\in S_c(\bO)$. For any $z\in \bO, N(z)\neq 0$ define
  $$g^z(u)=g\left(\frac{u z}{\N(z)}\right)|\N(z)|^{-4}.$$
Then $\mF_\bO(g^z)(x)=\mF_{\bO}(g)(xz)$. 
\end{Lem}

\begin{proof}
  $$\mF_{\bO}(g^z)(x)=
  \int\limits_{\bO} g^z(y)\psi( \T( y  \bar x)) dy=$$
  $$ \int\limits_{\bO} g\left(\frac{yz}{\N(z)}\right)
  \psi( \T\left(\frac{yz}{\N(z)}(\bar z \bar x)\right))
 |\N(z)|^{-4+4}  d\left(\frac{yz}{\N(z)}\right)=$$
 $$\int\limits_{\bO} g(y)\psi(\T(y\overline{xz})d y =\mF(g)(xz).$$
\end{proof}

 \section{Main Theorem}\label{sec:main}

\begin{Thm}\label{main}
  One has $\Pi(s)(f)=\Phi(f)$ for any $f\in \mS_c(\Omega)$.
\end{Thm}

The proof occupies the rest of the section.
We prove  Theorem \ref{main}, based on Propositions
\ref{Pi(s):messy} and \ref{Pi(s):decomposition} below.
The proof of Propositions follow the proof of Theorem \ref{main}.
\begin{proof}
\begin{enumerate}
\item In Proposition \ref{Pi(s):messy} we compute
  $\Pi(s)(f)=\Pi(s_1)\Pi(s_2)\Pi(s_3)(f)$
  for any $f\in \mS_c(\Omega^1),$
  using the formulas for $\Pi(s_0),\Pi(s_2)$ and $\Pi(s_3)$
  from section \ref{sec:minimal}.
  The goal is to bring the resulting expression to the form of $\Phi(f)$. 
  Using the decomposition in Theorem \ref{delta:decomp}
  for the quadratic space $(\bO, N)$, we deduce in 
  Corollary \ref{cor:T1:T2} the decomposition $\Pi(s)(f)=T_1(f)+T_2(f)$
  for all $f\in \mS_c(\Omega^1)$.
  
\item Proposition \ref{Pi(s):decomposition} consists of two parts. 
  In part $(2)$ we show that $T_2(f)=\Phi_2(f)$ for $f\in \mS_c(\Omega^1)$.
  This is the most involved computation in the proof, which uses
  at a final point the identity  from Lemma \ref{inner:product:identity}.  
It follows that  $\tilde T_1=\Pi(s)-\Phi_2$  is a $(M,s)$ equivariant
  operator on $\mS_c(\Omega)$, whose restriction to $\mS_c(\Omega^1)$
  equals to $T_1$.
\item Part  $(1)$ of Proposition \ref{Pi(s):decomposition} claims that
   $T_1(f)$ is contained in $\mS(\Omega,3)$. Since $\mS_c(\Omega^1)$
  generates $\mS_c(\Omega)$ as $M$-representation,
  it follows that $\tilde T_1,\Phi_1\in \Hom_{(M,s)}(\mS_c(\Omega),\mS(\Omega,3))$,
  that is one-dimensional by
  Proposition \ref{one:dim}. Hence
  $\tilde T_1= c\cdot\Phi_1$  for some $c\in \C$.
  Consequently
  \beq \label{Pi(s)-Phi}
  \Pi(s)-\Phi=(c-1)\cdot \Phi_1.
  \eeq
\item
  Plug the function $f_0$ from Proposition \ref{f0:lemma}
  into equation (\ref{Pi(s)-Phi}) and evaluate at a point
  $t\cdot c_0$ for $|t|$ large.
  Since $\Pi(s)(f_0)$ has bounded support, it follows that LHS is zero.
  On the other hand  $\Phi_1(f_0)(t c_0)=|t|^{-3}\Phi_1(f_0)(c_0)\neq 0,$
  which  implies that $c=1$.
 Hence $\Pi(s)=\Phi$ on $\mS_c(\Omega)$.
  This concludes the proof of Theorem \ref{main}.
  \end{enumerate}
\end{proof}
It remains to prove Propositions \ref{Pi(s):messy} and
\ref{Pi(s):decomposition}. The set $\Omega^1\subset \Omega$
is defined in (\ref{omega1}).

\begin{Prop}\label{Pi(s):messy}
  Let $c_1=c(a,x,\bar y)\in \Omega^1$ and
  $f\in \mS_c(\Omega_1)$.
  Then
  \beq
\Pi(s)(f)(c_1)=\int\limits_{F^\times} \left|\frac{N(x)t}{a}\right|^{-3} 
\left(\int\limits_{\bO(\frac{tN(x)}{a})} I(f)(x',t) \eta_{\frac{tN(x)}{a}}(x')\right) 
\cdot   |t|^{-4} d^\times t,
\eeq
where
\beq \label{I(f):def}
 I(f)(x',t)=\int\limits_\bO
\int\limits_\bO
f(c(t, u,\bar v))
\psi\left(\left\la \bar x',
\frac{x a}{N(x)}+
\frac{u}{t}+
\frac{\bar v(yx) }{N(x)t}\right\ra\right) du dv
\eeq
\end{Prop}

\begin{proof}
  We use the fact that $s=s_1s_2s_3=s_0s_2s_0s_2s_3$.
  The action of $s_0$ is given by
$$\Pi(s_0) f(c(a,x,\bar y))=
f\left(c(\frac{\N(x)}{a},\bar x , \frac{\bar x \bar y}{a})\right)$$
  and the action of
  $s_2,s_3$ can be read from Corollary   \ref{Schrodinger:action:s2:s3}.

 Hence  
  $$\Pi(s_1)(f)(c_1)=\Pi(s_0 s_2 s_0)(f)(c(a,x,\bar y))=
  \Pi(s_2 s_0)(f)\left(
  \frac{N(x)}{a}, \bar x , \frac{\bar x \bar y}{a}
  \right)=$$
  $$\left|\frac{N(x)}{a}\right|^{-4}
  \int\limits_{\bO}
  \Pi(s_0)(f) \left(\frac{N(x)}{a}, x',  \frac{\bar x \bar y}{a}\right)\,
  \psi\left(
  \frac{\la x', \bar x\ra a}{N(x)}
  \right) dx'=$$
 $$\left|\frac{N(x)}{a}\right|^{-4}
  \int\limits_{\bO}
  f\left(\frac{N(x')a}{N(x)}, \bar x',
  \frac{\bar x'(\bar x \bar y) }{N(x)}\right)\,
  \psi\left(\frac{\la \bar x',x\ra a}{N(x)}\right) dx'.$$ 
  Breaking the integral over $x'\in \bO$ over fibers of $N$
  we get
 $$\left|\frac{N(x)}{a}\right|^{-4} 
  \int\limits_F \int\limits_{\bO(t)}
  f\left(\frac{ta}{N(x)}, \bar x',
  \frac{\bar x'(\bar x \bar y) }{N(x)}\right)\,
  \psi\left(
  \frac{\la \bar x',x\ra a}{N(x)}\right) |\eta_t(x')|  dt=$$
Make a change of variables $\frac{ta}{N(x)}\mapsto t$.
$$\left|\frac{N(x)}{a}\right|^{-3} 
  \int\limits_F \int\limits_{\bO\left(\frac{tN(x)}{a}\right)}
  f\left(t, \bar x',  \frac{\bar x'(\bar x \bar y) }{N(x)}\right)\,
  \psi\left(
  \frac{\la \bar x',x\ra a}{N(x)}
  \right) |\eta_{\frac{tN(x)}{a}}(x')|  d t.$$   
Finally, we   compute $\Pi(s_1s_2s_3)(f)(c_1)$ for 
$c_1=c(a,x,\bar y)\in \Omega^1$. It equals
%
 $$\left|\frac{N(x)}{a}\right|^{-3}
  \int\limits_F \int\limits_{\bO\left(\frac{tN(x)}{a}\right)}
  \Pi(s_2 s_3)(f)
  \left(c(t, \bar x',  \frac{\bar x'(\bar x \bar y) }{N(x)})\right)\,
  \psi\left(
  \left\la \bar x',\frac{x a}{N(x)}\right\ra
  \right)
  \eta_{\frac{tN(x)}{a}}(x')  d t=$$ 
\begin{multline} \label{still:fine}
\int\limits_{F^\times} \left|\frac{N(x)t}{a}\right|^{-3} 
\int\limits_{\bO(\frac{tN(x)}{a})}\\
\int\limits_\bO
\int\limits_\bO
f(c(t, u,\bar v))
\psi\left(\left\la \bar x',
\frac{x a}{N(x)}+
\frac{u}{t}+
\frac{\bar v(yx) }{N(x)t}\right\ra\right)
du dv\\
 \eta_{\frac{tN(x)}{a}}(x')\cdot   |t|^{-4} d^\times t.
\end{multline}

After substituting the expressions for $\Pi(s_2),\Pi(s_3)$ we have used
     $$\la  \bar x'(\bar x\bar y), \bar v\ra =
   T(\bar x' (\bar x\bar y) v)=
 T(\bar x' \overline{\bar v (yx)})=\la\bar x', \bar v (yx)\ra.$$

The middle line of \ref{still:fine} is exactly $I(f)(x',t)$ and Theorem is proved. 
\end{proof}

We apply Theorem
\ref{delta:decomp} to
the quadratic space $(\bO,N)$,
or precisely the equation (\ref{O(x):decomp}) for $I(f)$ and 
 $s=\frac{tN(x)}{a}$ to conclude
\begin{multline} \label{I(f):decomposition}
\int\limits_{\bO(\frac{tN(x)}{a})} I(f)(x',t) \eta_{\frac{tN(x)}{a}}(x')=
\int\limits_{\bO(0)} I(f)(x',t) \eta_{0}(x')+\\
\left|\frac{tN(x)}{a}\right|^3
\int\limits_{\bO} \mF_{\bO}(I(f))(x',t)H_3\left(\frac{tN(x)N(x')}{a}\right) dx'.
\end{multline}

We denote the first and second summands by $I_1(f)(t)$ and $I_2(f)(t)$.

\begin{Cor}\label{cor:T1:T2} For any $f\in \mS_c(\Omega^1)$  there is a  decomposition $$\Pi(s)(f)=T_1(f)+T_2(f),$$
where
\begin{multline} \label{T_1:expression}
T_1(f)(c(a,x,\bar y))=
  \int\limits_{F^\times} \left|\frac{N(x)t}{a}\right|^{-3} 
\int\limits_{\bO(0)}\\
\int\limits_\bO
\int\limits_\bO
f(c(t, u,\bar v))
\psi\left(\left\la \bar x',
\frac{x a}{N(x)}+
\frac{u}{t}+
\frac{\bar v(yx) }{N(x)t}\right\ra\right)
du dv\\
 \eta_{0}(x')\cdot   |t|^{-4} d^\times t.
\end{multline}
and
\beq \label{T_2:expression}
  T_2(f)(c(a,x,\bar y))=
  \int\limits_{F^\times} \left( \left|\frac{N(x)t}{a}\right|^{-3} I_2(f)(t)\right)
  |t|^{-4} d^\times t.
  \eeq
 
  \end{Cor}

We remind the reader that
$$\Phi_2(f)(c_1)=\int\limits_{F^\times}
\left(
\int\limits_{F^\times}\int\limits_\bO\int\limits_\bO
f(c_2)\psi(\la c_1,c_2\ra_J\cdot r) dudv (\psi(r^{-1})-1)|r|^3 d^\times r
\right)
|t|^{-4} d^\times t.$$

\begin{Prop}\label{Pi(s):decomposition}
  Let  $f\in \mS_c(\Omega^1)$.
  \begin{enumerate}
     \item $T_1(f)\in \mS(\Omega,3),$  
  \item    $T_2(f)=\Phi_2(f)$.
\end{enumerate}     
\end{Prop}

\begin{proof}
It is easy to see from (\ref{T_1:expression}) that
$T_1(f)(c(ba,bx,b\bar y))=|b|^{-3}T_1(f)(c(a,x,\bar y))$
for any $b\in F^\times$.
Equivalently  $\lambda(b) T_1(f)=|b|^3 T_1(f)$  and hence
$T_1(f)\in \mS(\Omega,3)$. 

To verify that $T_2(f)(c_1)=\Phi_2(f)(c_1)$ for $c_1=c(a,x,\bar y)$ it is enough to show that the summand $I_2(f)(t)$ in \ref{I(f):decomposition}.
 equals
\beq \label{I(f):second:summand}
\left|\frac{tN(x)}{a}\right|^3
\int\limits_{F^\times}\int\limits_{\bO}\int\limits_{\bO}f(c_2)
  \psi(\la c_1,c_2\ra \cdot r) du dv\cdot \left(\psi(r^{-1})-1\right)|r|^3 d^\times r,
  \eeq
where $c_2=c(t,u,\bar v)$.

Let us prove (\ref{I(f):second:summand}). Let $f=f_1\otimes f_2\otimes f_3$, where $f_1\in \mS_c(F^\times)$ and
   $f_2,f_3\in \mS_c(\bO^\times)$.
   In particular $f(c(t,u,\bar v))=f_1(t)f_2(u)f_3(\bar v)$. 

Rewriting the expression from (\ref{I(f):def}) we see that $I(f)(x',t)$ equals 
$$f_1(t)\cdot
\int\limits_\bO f_2(u) \psi\left(\frac{\la \bar x',u\ra }{t}\right) du\cdot 
\int\limits_\bO f_3(\bar v)
\psi\left(\frac{\la \bar x'(\bar x\bar y), \bar v\ra }{N(x)t}\right) d\bar v
\cdot 
\psi\left(\frac{ \la\bar x',x\ra a}{N(x)}\right)=$$

$$f_1(t)\cdot \mF_{\bO}(f_2)\left(\frac{\bar x'}{t}\right)\mF_{\bO}(f_3)
\left(\bar x'\frac{\bar x\bar y} {N(x)t}\right)\cdot
\psi\left(\frac{\la\bar x',x\ra a}{N(x)}\right).$$

Define   
$$\tilde f_2(u)=f_2(tu)|t|^8, \quad
\tilde f_3(u)=f_3\left(u \frac{\bar x \bar y t}{N(y)}\right)
\left|\frac{N(y)}{N(x)t^2}\right|^{-4}$$
Then by Lemma \ref{Fourier:O:shift} one has 
  $$\mF_{\bO}(\tilde f_2)(\bar x')=\mF_{\bO}(f_2)\left(\frac{\bar x'}{t}\right),\quad
\mF_\bO(\tilde f_3)(\bar x') =
\mF_\bO(f_3)\left(\bar x'\frac{\bar x \bar y}{N(x)t}\right) .$$
In addition,
$\mF_\bO(\tilde f_2)(\bar x') \mF_\bO(\tilde f_3)(\bar x')=
\mF_\bO(\tilde f_2\ast \tilde f_3)(\bar x').$
Hence 
$$I(f)(x',t)=f_1(t)\cdot  \mF_\bO(\tilde f_2\ast \tilde f_3)(\bar x')\cdot
\psi\left(\frac{ \la\bar x',x\ra a}{N(x)}\right).$$

Consequently, using equation \ref{O(x):decomp} for
$s=\frac{tN(x)}{a}$ and $I(f)$ we have

\begin{multline}
 I_2(f)(t)=\left|\frac{N(x)t}{a}\right|^{3}\cdot  f_1(t)
\int\limits_{\bO} (\tilde f_2\ast \tilde f_3)(u)
H_3\left(\frac{t\N(x)}{a} \N\left(u+\frac{ ax }{\N(x)}\right)\right)
du.
\end{multline}
After substituting the definition of $H_3$, and the definition of the convolution this becomes
\begin{multline}\label{last}
  \left|\frac{N(x)t}{a}\right|^{3}\cdot  \int\limits_{r\in F^\times}\\
  \int\limits_{\bO} \int\limits_{ \bO} f_1(t)\tilde f_2(u-\bar v) \tilde f_3(\bar v)
  \psi\left(\frac{t\N(x)}{a} \N\left(u+\frac{ ax }{\N(x)}\right) r\right)
  du dv \\
   (\psi(r^{-1})-1) |r|^3 d^\times r.
\end{multline}

It remains to show that the second line of this equation equals 
$$\int\limits_{\bO}\int\limits_{\bO}f(c(t,u,\bar v))
  \psi(\la c_1,c_2\ra_J r) du dv.$$

Indeed, the second line equals
$$\int\limits_{\bO} \int\limits_{\bO}
f_1(t)\tilde f_2(u) \tilde f_3(\bar v)
\psi\left(\frac{t\N(x)}{a}
\N\left(u+\bar v+\frac{ax }{\N(x)}\right)\cdot r\right) du dv =$$
$$\int\limits_{\bO} \int\limits_{\bO}
f_1(t)f_2(tu)  f_3\left(v\frac{\bar x \bar y t}{N(y)}\right)
\psi\left(\frac{tN(x)}{a}
\N\left(u+\bar v+\frac{ax }{\N(x)}\right)\cdot r\right)
|t|^8  du \left|\frac{\N(y)}{\N(x)t^2}\right|^{-4} dv =$$

$$\int\limits_{\bO} \int\limits_{\bO}
f_1(t)f_2(u)f_3(\bar v)
\psi\left(\frac{t\N(x)}{a}
\N\left(\frac{u}{t}+\frac{\bar v (yx)}{t\N(x)} +\frac{ax }{\N(x)}\right)\cdot r\right)  du dv=$$
$$\int\limits_{\bO} \int\limits_{\bO}
f(c(t,u,\bar v))
\psi\left(\frac{t\N(x)}{a}
\N\left(\frac{u}{t}+\frac{\bar v (yx)}{t\N(x)} +\frac{ax }{\N(x)}\right)\cdot r\right)  du dv.$$
We are done after we prove  the following remarkable identity.

\begin{Lem}\label{inner:product:identity}
  Let   $c_1=c(a,x,\bar y)$ and $c_2=c(t,u,\bar v)$
  be two elements in $\Omega^1.$
  $$ \frac{t\N(x)}{a}\cdot \N\left(\frac{ax}{N(x)}+\frac{u}{t}+\frac{\bar v(yx) }{N(x)t}\right)
  = \la c_1,c_2\ra_J.$$
  \end{Lem}

 \begin{proof}  
 One has 
  $$\frac{ax}{N(x)}+\frac{u}{t}+\frac{\bar v(yx)}{N(x)t}=
  \frac{(at+u\bar x+((\bar v(yx))x^{-1}))x}{N(x)t}.$$
  Using the multiplicavity of the norm it is enough to show that 

  \beq \label{modified:identity}
(at)^{-1}\cdot N(at+u\bar x+ (\bar v(yx))x^{-1})=\la c_1,c_2\ra_J.
\eeq
To prove \ref{modified:identity} we  write explicitly the right hand side.
By Lemma \ref{pairing:trace} one has 
  $$\la c_1,c_2\ra_J=\bT(c_1\circ c_2)=$$
  $$at +\frac{N(x)N(u)}{at}+ \frac{N(y)N(v)}{at}+
  T(\bar x u)+T(\bar v y)+ \frac{T((yx)(\bar u \bar v))}{at}.$$

  Next we elaborate the left-hand side. Let
  $w=at+u\bar x+ (\bar v(yx))x^{-1}$. Then
  $$\frac{N(w)}{at}= \frac{1}{at}
\left(at+u\bar x+ (\bar v(yx))x^{-1}\right)\cdot
\left(at+x\bar u+ \bar x^{-1}((\bar x\bar y)v)\right)
=$$
  $$at+\frac{N(x)N(u)}{at}+ \frac{N(y)N(v)}{at} +T(\bar x u)+
   T(\bar v(yx)x^{-1})+ \frac{T(((\bar v(yx))x^{-1})(x\bar u)) }{at}.$$
  
The first four terms for RHS and LHS are the same. 
Any subalgebra generated by two elements in $\bO$ is associative
and in particular $(yx)x^{-1}=y$.  This implies
$\T(\bar v(yx)x^{-1})=\T(\bar v((yx)x^{-1}))=T(\bar vy)$.
For the same reason
$$\T(((\bar v(yx))x^{-1})x\bar u)=
\T(\bar u(((\bar v(yx))x^{-1}) x))=\T(\bar u \bar v (yx))=
\T\left((yx)(\bar u\bar v)\right)$$
and so all the six terms on both sides are equal.
\end{proof} 

\end{proof}

The identity \ref{modified:identity} in $H_3(\bO)$
is quite subtle, since  the geometric meaning
of the LHS is not clear to us. We note that the analogous identity exists
for elements of rank one in the Jordan algebra $H_3(\bB)$, where 
$\bB$ is an {\sl associative} composition algebra. In this case there is 
a simple explanation for the identity. So assume that 
the elements $x,y,u,v$ belong to $\bB$. The identity (\ref{modified:identity}) now reads
\beq \label{modified:identity:B}
(at)^{-1}\cdot \N(at+u\bar x+ \bar v y)=\la c_1,c_2\ra_J.
\eeq

The elements $c(a,x,\bar y), c(t,u,\bar v)$ can be written as 
$$c(a,x,\bar y)=a\cdot \bar \bv^t\cdot \bv, 
\quad c(t,u,\bar v)=t\cdot \bar \bw^t\cdot \bw$$
where 
$$\bv=(1, \frac{x}{a},\frac{\bar y}{a}), \quad \bw=(1, \frac{u}{t},\frac{\bar v}{t})\in \bB^3, \quad a,t\in F^\times .$$
The left hand side of \ref{modified:identity:B} is $at\cdot N(\bv \bar \bw^t)$.

From this 
\begin{multline}
    \bT(c_1\circ c_2)= \frac{1}{2} at\cdot  
\bT((\bar \bv^t \bv)(\bar\bw^t  \bw)+(\bar \bw^t \bw)(\bar \bv^t\bv))=\\
\frac{1}{2}at\cdot 
\bT(\bar \bv^t (\bv\bar\bw^t)  \bw)+(\bar \bw^t (\bw \bar \bv^t)\bv))=at\cdot\N(\bv\bar \bw^t)
\end{multline}
as required.

\section{The Jacquet Functor of $\Pi$}\label{sec:jacquet}

It is known, see \cite{Sa2}, that the Jacquet module $\Pi_N$ as a representation
of $M$ is {\sl isomorphic} to a direct sum of a
one-dimensional representation and a representation of $M$,
whose restriction to $[M,M]$ is the unique minimal representation
$\Pi_1$ of $E_6$. This fact is often used for establishing
the theta correspondence for various dual pairs in $H$. 

Using the realization $\Pi$ in $\mS(\Omega)$ it is easy to see
that $\Pi_N$ is canonically identified with the space of germs
$[\mS(\Omega)]_0$ at zero. As an application of Theorem \ref{main} we can show that
$[\mS(\Omega)]_0$ naturally breaks in to this  direct sum.

\begin{Prop}
  $\mS(\Omega)=\mS_c(\Omega)+\Pi(s)\mS_c(\Omega).$
\end{Prop}
\begin{proof}
  It is clear that the RHS is contained in LHS and  is $\Pi(s)$ invariant.
  For any vector $f$ in RHS
  and $n\in N$, one has $\Pi(n)f-f\in \mS_c(\Omega)$. Hence RHS is
  $\Pi(n)$ invariant. Since $N,s$ generate $H$, it follows that RHS is
  $H$-invariant subspace of LHS
  and hence equals to the irreducible representation $        \mS(\Omega)$. 
  \end{proof}

Hence any element in $[\mS(\Omega)]_0$ has a representative $\Phi(f)$ for
$f\in \mS_c(\Omega)$. Recall that $\Phi_2(f)$ is constant in a neighborhood
of zero. We denote this constant by $\Phi_2(f)(0)$. 

\begin{Cor}
The map
$$[\mS(\Omega)]_0\rightarrow  |\omega_7|^4 \oplus |\omega_7|^{-6}\otimes \Pi_1$$
defined by
$$\Phi(f)\mapsto (\Phi_2(f)(0),\Phi_1(f))$$
for any $f\in \mS_c(\Omega)$ is well-defined $M$-isomorphism. 
\end{Cor}

\begin{proof} The map is obviously $M$-equivariant. It is injective. Indeed, if $\Phi_2(f)(0)=0$ and $\Phi_1(f)=0,$ then $\Phi(f)$
vanishes in a neighborhood of zero, i.e. represents a zero germ. It is surjective  since $\Phi_1,\Phi_2(\cdot)(0) $ are non-zero maps 
 and the codomain is a sum of two distinct irreducible representations.
\end{proof}

\section{Braverman-Kazhdan operators}\label{sec:BK}
We start with an overview of an important work by Braverman and Kazhdan.
Our notation is compatible with \cite{BravermanKazhdan02}
and is not related to our former notation.

Let $G$ be the group of $F$-points of
a split simply-connected group, $M$ be a Levi subgroup of $G$
and $M^{ab}=M/[M,M]$.
For any parabolic subgroup $P$ containing $M$ as a Levi factor,
we may consider $X_P=[P,P]\backslash G$, that carries a natural
action of $G\times M^{ab}$. There is a  $G$-invariant measure on
$X_P$, unique up to constant.  The space $L^2(X_P)$ carries a unitary
action of $G\times M^{ab}$ via
$$(g,m)f(x)=f(mxg)|\delta_P(m)|^{1/2}.$$

The main result of \cite{BravermanKazhdan02}
is a construction of unitary $G\times M^{ab}$ equivariant operators
$\mF_{P,Q}:L^2(X_P)\rightarrow L^2(X_Q)$ for any two parabolic
subgroups $P,Q$ sharing the Levi factor $M$. It is shown that
$\mF_{P,P}=Id$ and $\mF_{Q,R}\circ \mF_{P,Q}=\mF_{P,R}$.
In particular $\mF_{Q,P}\circ \mF_{P,Q}=Id$.
These operators, that are in fact normalized intertwining operators,
are also called generalized Fourier transforms.
A motivating  example is the case $G=SL(V)$,
for a finite-dimensional vector space $V$,
and $P$,$Q$ are  maximal parabolic subgroups  stabilizing a line in $V$
and a line in $V^\ast$ in general position.
Then $X_P=V-0$ and $X_Q=V^\ast-0$. It is easy to check that
the construction in \cite{BravermanKazhdan02}
gives the classical Fourier transform.
$$\mF_{P,Q}:L^2(V)\rightarrow L^2(V^\ast),\quad
\mF_{P,Q}(f)(v^\ast)=\int\limits_V f(v)\psi(\la v^\ast,v\ra) dv.$$ 
Let us give some details on the construction.
One begins by defining the Radon transforms
$$\mR_{P,Q}: \mS_c(X_P)\rightarrow \mS^\infty(X_Q),$$
that are $G\times M^{ab}$ equivariant, but not unitary.
They can be thought as unnormalized intertwining operators.
$$\mR_{P,Q}(f)(x)=\int\limits_{(x,y)\in Z_{P,Q}} f(x) dx$$
where $Z_{P,Q}\subset X_P\times X_Q$ is the image of the projection
$G\rightarrow X_P\times X_Q$.
Next one defines a distribution $\eta_{P,Q}$ on the torus
$M^{ab}$ with Langlands dual $Z(M^\vee)$, and puts
$$\mF_{P,Q}=\eta_{P,Q}(\mR_{P,Q}(f)).$$
We refer the reader to \cite{BravermanKazhdan02} for the definition of $\eta_{P,Q}$
for general $P,Q$. Below we apply the recipe for our particular case. 

To relate the construction to our work we shall describe
the operator $\mF_{Q,\bar Q}$, where $G$ is a a split group
of type $E_6$ and $Q=MU$ is a maximal
parabolic subgroup of type $D_5$. The group $M^{ab}$ is a one-dimensional
torus, generated by $\omega_6^\vee$.
In this case $X_Q$ is identified with $\Omega\subset \bar{\mf{n}}$ and
$X_{\bar Q}$ is identified with $\bar \Omega\subset \mf{n}$, the
elements of rank $1$ in the associated Jordan algebra. 
The map $\kappa:\Omega\times \bar\Omega\rightarrow F$
defines a $G$-invariant pairing, so that
$$Z_{Q,Q'}=\{(x,y): \kappa(x,y)=1\}\subset X_Q\times X_{\bar Q}.$$

For any $y\in \bar \Omega$  consider a fibration
$$p_y:\Omega\rightarrow F, \quad p_y(x)=\kappa(x,y)$$
and denote by $\Omega_a(y)$ the fiber over $a$ with the induced
measure $\eta_{y,a}$. 
  The Radon transform in this case is  
  $$\mR_{Q,\bar Q}(f)(y)=\int\limits_{\Omega_1(y)} f(x)\eta_{1,\eta}(x)$$
  and
  $$\mR_{Q,\bar Q}(a f)(y)=
  \int\limits_{\Omega_a(y)} f(x)\eta_{a,\eta}(x).$$
  

  Let us apply the recipe in \cite{BravermanKazhdan02}
  for the distribution $\eta_{Q,\bar Q}$.
The space ${\mf{u}_{Q,\bar Q}^\vee=\mf{u}^\vee}$ has dimension $16$.
The principal subalgebra $\mf{sl}_2$ in $\mf{m}$
acts on $\mf{u}^\vee$, which is a sum of two irreducible
representations of dimension $n_1=11$ and $n_2=5$.
The one-dimensional torus  $Z(M^\vee)=(M^{ab})^\vee$
acts on the two-dimensional space of highest vectors
in $\mf{u}^\vee$ with characters
$$\lambda_1, \lambda_2\in \Hom(Z(M^\vee), \bG_m)=\Hom(\bG_m,M^{ab}).$$

The distribution $\eta_{Q,\bar Q}$ on $M^{ab}$ is defined as a convolution
\beq 
(\lambda_1)_{!}(\psi(t)|t|^{\frac{n_1}{2}})\ast
(\lambda_2)_{!}(\psi(t)|t|^{\frac{n_2}{2}}).
\eeq
Equivalently, as a distribution on $M^{ab}\simeq F^\times$
$$\eta_{Q,\bar Q}(a)=\int\limits_{F^{\times}} \psi(ar)\psi(r^{-1})|r|^{\frac{11-5}{2}}
d^\times r$$
Hence
$$\mF_{Q,\bar Q}(f)(y)=
\int\limits_{M^{ab}} \mR_{Q,\bar Q}(af)\cdot \eta_{Q,\bar Q}(a) d^\times a=$$
$$\int\limits_{F^\times}
\int\limits_{\Omega_{a(y)}} f(x) \eta_{y,a}(x)\cdot 
\psi(ar)\psi(r^{-1}) |r|^{3}d^\times r=$$
$$\int\limits_F
\left(\int\limits_\Omega f(x)\psi(\kappa(x,y) \cdot r)\eta(x)\right)
\psi(r^{-1})|r|^3 d^\times r.$$

Note, that while the $\mF_{Q,\bar Q}(f)(y)$ converges only for $f$ in some
dense subset in $L^2(X_Q)$, the 
last expression makes sense for all $f\in \mS_c(\Omega)$. 

Recall  that we have defined a pairing on $\Omega$ by 
$\la x, y\ra=-\kappa(x,sys^{-1})$.

Hence 
$$\mF_{Q,\bar Q}(f)(sys^{-1})=\int\limits_F
\left(\int\limits_\Omega f(x)\psi(-\la x,y\ra \cdot r)\eta(x)\right)
\psi(r^{-1})|r|^3 d^\times r$$
which is exactly the expression for $\Pi(s\lambda(-1))f(y)$. 

To sum up, denoting by $\iota_s$ the unitary operator
$$\iota_s:L^2(X_{\bar Q})\rightarrow L^2(X_Q),
\quad \iota_s(f)(y)=f(sys^{-1})$$
and $\tilde s=s\lambda(-1)$,
we have the relation
\begin{Thm} 
  $\iota_s\circ \mF_{Q,\bar Q}=\Pi(\tilde s)$.
\end{Thm}

It is remarkable that the space $L^2(X_Q)$ besides
being a natural representation of $G\times M^{ab}$
can be extended to the minimal representation of a bigger group $H$
and a distinguished involution $s$ in it acts by a (twisted)
generalized Fourier transform. This phenomenon can be observed
in the example of the minimal representation
of the double cover of the symplectic group, $G=Sp(V\oplus V^\ast)$
with $M=GL(V)$, realized on $L^2(V)$ 
where a representative of the longest element in the Weyl group acts by the classical
Fourier transform. Similarly, for a group of type $D_n$,
the minimal representation is realized in a cone  of isotropic vectors
in a space of dimension $2n-2$. The formula for the operator
$\Pi(s)$, called a Fourier transform on a cone, has been obtained
in \cite{GurevichKazhdan}. In a forthcoming work, we shall treat
the case  $G=E_6$. 
\section*{Acknowledgement}
 This project was initiated when the authors participated in a Research in Teams (RIT) project 
“Modular Forms and Theta Correspondence
for Exceptional Groups” hosted by the Erwin Schrodinger Institute, University of Vienna.
We thank the Erwin Schrodinger Institute for its generous support and wonderful working
environment, as well as our other RIT members  Aaron Pollack 
and Gordan Savin for helpful discussions. We thank David 
Kazhdan for inspiring ideas and stimulating conversation.
The second author 
gratefully acknowledges National  University of Singapore for its warm hospitality and the excellent working conditions provided during the completion of this project. 
W.T. Gan is supported by a Tan Chin Tuan Centennial Professorship at the National University of Singapore.
N. Gurevich is partially supported by ISF grant 1643/23.

\bibliographystyle{alpha}
 \bibliography{bib}

\end{document}